\documentclass[a4paper,10pt]{siamltex}
\usepackage{cite}
\usepackage{amstext}
\usepackage{amssymb}
\usepackage{amsmath}
\usepackage{appendix}
\usepackage{amscd}
\usepackage{subfloat}
\usepackage{subcaption}
\usepackage[dvips]{graphicx}
\usepackage[latin1]{inputenc}
\usepackage[english]{babel}

\usepackage{graphicx}
\usepackage{tabularx}
\usepackage{multirow}
\usepackage{bbm} 
\usepackage{todonotes}

\usepackage{algorithm}
\usepackage{verbatim}
\usepackage{color}
\usepackage{xcolor}
\usepackage{geometry}

\usepackage{enumerate}
\usepackage{ulem}

\newtheorem{assumption}[theorem]{Assumption}

\usepackage{hyperref}

\newcommand{\be}{\begin{equation}}
\newcommand{\ee}{\end{equation}}
\newcommand{\ba}{\begin{array}}
\newcommand{\ea}{\end{array}}
\newcommand{\bee}{\begin{eqnarray*}}
\newcommand{\eee}{\end{eqnarray*}}
\newcommand{\bea}{\begin{eqnarray}}
\newcommand{\eea}{\end{eqnarray}}
\newcommand{\eqdef}{\stackrel{\rm def}{=}}

\newcommand{\R}{\mathbb{R}}
\newcommand{\E}{\mathbb{E}}

\newcounter{algo}[section]
\renewcommand{\thealgo}{\thesection.\arabic{algo}}
\newcommand{\algo}[3]{\refstepcounter{algo}
\begin{center}
\framebox[\textwidth]{
\parbox{0.95\textwidth} {\vspace{\topsep}
{\bf Algorithm \thealgo . #2}\label{#1}\\
\vspace*{-\topsep} \mbox{ }\\
{#3} \vspace{\topsep} }}
\end{center}}
\def\IR{\hbox{\rm I\kern-.2em\hbox{\rm R}}}
\def\IC{\hbox{\rm C\kern-.58em{\raise.53ex\hbox{$\scriptscriptstyle|$}}
    \kern-.55em{\raise.53ex\hbox{$\scriptscriptstyle|$}} }}
    
\def\IR{\hbox{\rm I\kern-.2em\hbox{\rm R}}}
\def\IC{\hbox{\rm C\kern-.58em{\raise.53ex\hbox{$\scriptscriptstyle|$}}
    \kern-.55em{\raise.53ex\hbox{$\scriptscriptstyle|$}} }}

\newcommand{\NN}{\mathbb{N}}
\newcommand{\RR}{\mathbb{R}}

\def \P {{\mathbb P}}
\def \E {{\mathbb E}}

\usepackage{bbm}				
\def \1{{\mathbbm 1}}

\newcommand{\hathat}[1]{\widehat{#1}^R}

\DeclareMathOperator{\spa}{span}
\DeclareMathOperator{\range}{range}

\DeclareMathOperator{\eig}{eig}
\DeclareMathOperator{\dist}{dist}

\title{A variable dimension sketching strategy for nonlinear least-squares}
 \author{ Stefania Bellavia\footnotemark[1], Greta Malaspina\footnotemark[1], Benedetta Morini\footnotemark[1]}

\begin {document}
\maketitle
\footnotetext[1]{Dipartimento  di Ingegneria Industriale, Universit\`a degli Studi di Firenze,
Viale G.B. Morgagni 40,  50134 Firenze,  Italia. Members of the INdAM Research Group GNCS. Emails:
stefania.bellavia@unifi.it, greta.malaspina@unifi.it,  benedetta.morini@unifi.it}
\footnotetext[2]{ 
The research that led to the present paper was partially supported by INDAM-GNCS through Progetti di Ricerca 2023 and 
by PNRR - Missione 4 Istruzione e Ricerca - Componente C2 Investimento 1.1, Fondo per il Programma Nazionale di Ricerca e Progetti di Rilevante Interesse Nazionale (PRIN) funded by the European Commission under the NextGeneration EU programme, project ``Advanced optimization METhods for automated central veIn Sign detection in multiple sclerosis from magneTic resonAnce imaging (AMETISTA)'',  code: P2022J9SNP,
MUR D.D. financing decree n. 1379 of 1st September 2023 (CUP E53D23017980001), project 
``Numerical Optimization with Adaptive Accuracy and Applications to Machine Learning'',  code: 2022N3ZNAX
     MUR D.D. financing decree n. 973 of 30th June 2023 (CUP B53D23012670006), and by Partenariato esteso FAIR ``Future Artificial Intelligence Research'' SPOKE 1 Human-Centered AI. Obiettivo 4, Project ``Mathematical and Physical approaches to innovative Machine Learning technologies (MaPLe)'', Codice Identificativo EP\_FAIR\_002, CUP B93C23001750006.}

\begin{abstract}
We present a stochastic inexact Gauss-Newton method for the solution of nonlinear least-squares. To reduce the computational cost with respect to the classical method, at each iteration the proposed algorithm approximately minimizes the local model on a random subspace. The dimension of the subspace varies along the iterations, and two strategies are considered for its update: the first is based solely on the Armijo condition, the latter is based on information from the true Gauss-Newton model. Under suitable assumptions on the objective function and the random subspace, we prove a probabilistic bound on the number of iterations needed to drive the norm of the gradient below any given threshold. Moreover, we provide a theoretical analysis of the local behavior of the method. The numerical experiments demonstrate the effectiveness of the proposed method.
\end{abstract}
\section{Introduction}
In this paper, we consider second-order methods for solving the nonlinear least-squares problem
\be \label{merit}
\min_{x\in\RR^n} f(x)=\frac 1  2 \|F(x)\|_2^2,
\ee
where the residual function  $F\colon \RR^n\rightarrow\RR^m$ is continuously differentiable, and the objective function $f\colon \RR^n\rightarrow\RR$  has a large number $n$ of variables. Our focus is on using second-order models of reduced dimension which still incorporate some form of curvature information. To this end, we propose a sketching strategy of variable dimension.

In recent years, randomized linear algebra \cite{Mahoney_randalg, acta, tropp, Wood} has emerged  as a powerful tool for solving   optimization problems with high computational and memory demands. Randomized sampling and randomized embeddings are the core of a variety of optimization methods with stochastic models that are suitable for solving many applications, including machine learning; see e.g.,  \cite{bgmt4, bmm, sirtr_coap, Berahas, STORM2, noc1, noc2, noc3, cfs, cfs2, CST, cartis, STORM1, PS, Pilanci, Shao21, sketched}. Referring to our problem (\ref{merit}), a large reduction of either the variable dimension $n$ or the dimension $m$ of the observations can be achieved via randomized linear algebra. The application of a random embedding, referred to as sketching, can be used to restrict the computation of the trial step to a random 
subspace of $\R^n$ of dimension considerably smaller than $n$. As a result, the per-iteration cost is reduced and savings occur in terms of both cost and memory. We refer to \cite[\S 1]{CST} for a detailed overview of the existing literature of optimization methods based on sketching.

We turned our attention to sketching methods motivated by the works \cite{cfs, cfs2, Shao21} where a general random subspace framework for unconstrained nonconvex optimization was proposed and then specialized  to trust-region and quadratic regularization methods applied to problem (\ref{merit}). Under the assumption that the random subspace is an embedding of the gradient of $f$ at the current iteration, such methods drive the gradient of $f$ below a specified threshold $\epsilon$, with high probability, in a number of iterations of order ${\cal{O}}(\epsilon^{-2})$. While sketching preserves the order of worst-case complexity, the numerical experience reported in \cite{cfs} is not conclusive on the benefit of using sketching. 

The need for further investigation on the performance of randomized subspace methods led us to develop a variable dimension sketching strategy for  Levenberg-Marquardt type models. Our main
effort is to analyze if the trial step captures second-order information though it belongs to the random subspace generated by sketching. To this end we need to assume that the sketching matrices embed the transpose of the Jacobian of $F$ at the iterates with some probability; consequently, we focus on cases where the Jacobian is low-rank, e.g., strongly underdetermined problems.
Our algorithm employs the Armijo condition to test the acceptance of the trial step and to adapt the size of the random subspace; it  belongs to the class of step search methods since the step direction can change during the back-tracking procedure \cite{JSX} and retains  
worst-case iteration complexity of order ${\cal{O}}(\epsilon^{-2})$. The choice of the size of the random subspace is made in view of two conditions; the first condition is the acceptance or rejection of the iterate based on the Armijo condition, the latter condition checks if the trial step captures second-order information. The numerical experience shows that our strategy improves the performance of the algorithm with respect to the basic application of sketching. We remark that we are aware of only one recent work concerning variable size strategies. In \cite{CST}, a cubic regularization method for unconstrained optimization problems is proposed along with variable size sketching matrices; the approach adopted  in \cite{CST} is very different from ours as
the procedure for choosing the size of the random matrices is based on the rank of the sketched Hessian matrix.
This work is organized as follows. In Section \ref{sec2} we present our algorithm and in Section \ref{sec3} we study its theoretical properties. In Section \ref{sec4} we discuss the effect of sketching on the approximate minimization of the full Levenberg-Marquardt model and introduce a practical strategy for choosing the random embedding size adaptively. In Section \ref{sec5} we analyze the local convergence behavior of a variant of our algorithm to further support our proposal for the adaptive choice of the embedding size.  Finally, in Section \ref{sec6} we present numerical results that illustrate the benefit of using our variable dimension sketching strategy. The appendix summarizes some matrix distributions from the literature that are of interest for our purposes and contains some proofs.

\section{A step search  algorithm with random reduced models}\label{sec2}
We introduce our Algorithm \ref{algo_general} 
 employing reduced models generated by randomized embedding and a step search strategy.
Given the iterate $x_k\in \R^n$, the reduced model is built relying on the Gauss-Newton model 
\begin{equation}\label{mk}
{m}_k(s)= \frac 1 2 \|J(x_k)s+F(x_k)\|_2^2,
\end{equation}
where $s\in \R^n$.  Suppose that a sketching matrix, i.e., a random  matrix $M_k\in \R^{\ell_k\times n}$, is sampled from a matrix distribution $\mathcal{M}_k$  and let  $s_k=M_k^T  \widehat s_k$, $\widehat s_k\in \R^{\ell_k}$. Then we can form the randomly generated model
\begin{equation}\label{mhat}
\widehat{m}_k(\widehat s)= \frac 1 2 \|J(x_k)M_k^T \widehat s+F(x_k)\|_2^2,
\end{equation}
whose dimension
is reduced whenever $\ell_k<n$. In what follows  we assume that $ \ell_k \in [\ell_{\min}, \ell_{\max}]$, $\ell_{\max}\le n$ and   we make use of  the  Levenberg-Marquardt model  of  the form
\begin{eqnarray}
 \min_{\widehat s\in \R^{\ell_k}}  \hathat{m}_k(\widehat s)&=&\frac 1 2 \|J(x_k)M_k^T \widehat s+F(x_k)\|_2^2+\frac 1 2 \mu_k\|\widehat s\|_2^2   \label{modelGN_S}  \\
&=& \frac 1 2 \left\|
\left(
\begin{array}{c} J(x_k)M_k^T \\ \sqrt{\mu_k} I_{\ell_k}  \end{array} 
\right)\widehat s+
\left(
\begin{array}{c} F(x_k) \\ 0 \end{array} 
\right)   \right\|_2^2, \nonumber
\end{eqnarray}
with $0<\mu_{\min}\leq\mu_k\leq \mu_{\max}$. Note that $\hathat m_k$ is a model for (\ref{merit})
in the subspace generated by $M_k^T$ and is strictly convex due to the regularization term. 
We observe that
\begin{eqnarray*}
  \nabla \widehat m_k(\widehat s)   &=& M_k J(x_k)^T J(x_k)M_k^T\widehat s_k+M_k J(x_k)^T F(x_k)\\
  \nabla \hathat m_k(\widehat s)   &=& \left(M_k J(x_k)^T J(x_k)M_k^T+\mu_k I_{\ell_k}\right) \widehat s+M_k J(x_k)^T F(x_k)
\end{eqnarray*}
and the
optimality condition $\nabla \hathat m_k(\widehat s)=0$ for (\ref{modelGN_S})  amounts to solving the linear system of dimension  $\ell_k\times \ell_k$ 
\begin{equation}\label{nw_sketch}
\left(M_k J(x_k)^T J(x_k)M_k^T+\mu_k I_{\ell_k}\right) \widehat s =-M_k J(x_k)^T F(x_k).
\end{equation}
which can be solved approximately finding a step $\widehat s_k$ s.t.
\begin{eqnarray}
& &   r_k=\left(
\begin{array}{c} J(x_k)M_k^T \\ \sqrt{\mu_k} I_{\ell_k} \end{array}\right)  \widehat s_k  + \left(
\begin{array}{c} F(x_k) \\ 0 \end{array} 
\right) ,
 \label{residuo}\\
& & (M_k J(x_k)^T J(x_k)M_k^T+\mu_k I_{\ell_k}) \widehat s_k =-M_k J(x_k)^T F(x_k)+\rho_k,\label{newin}\\
& &  \|\rho_k\|_2=\|(M_k J(x_k)^T \  \sqrt{\mu_k} I_{\ell_k} ) r_k \|_2\le \eta_k \|M_k J(x_k)^T F(x_k)\|_2,\label{stepin}
\end{eqnarray}
for some $0\le \eta_k\leq\eta_{\max}<1$.  Choosing $\eta_k=0$ corresponds to finding the exact minimizer of $\hathat m_k$, otherwise the step is an approximate minimizer computed using an iterative solver, and in this work we adopt  LSMR \cite{LSMR}. Once $\widehat s_k$ is available, the trial step $s_k=M_k^T \widehat s_k$ in the full space is recovered.
The procedure above describes Step 1 of Algorithm \ref{algo_general}.
\vskip 10pt
\algo{}{General scheme: $k$-th iteration }{
 \noindent
\hspace*{1pt} Given  $c, \gamma,\widehat\gamma \in (0,1)$, $t_{\max}, \eta_{\max}, {\mu}_{\min}, \mu_{\max}>0$, $\ell_{\min}, \ell_{\max}\in \mathbb{N}$,
$\ell_{\min}<\ell_{\max}\le n$.
\\
\hspace*{1pt} Given   $x_k\in {\mathbb R}^n$, $t_k \in(0, t_{\max}]$, $\eta_k\in [0,  \eta_{\max}]$,   $\mu_k\in [{\mu}_{\min}, \mu_{\max}]$,
$\ell_k \in \mathbb{N},\,  \ell_k\in [\ell_{\min},\ell_{\max}]$.
\vskip 2 pt
 \begin{description}
 \item{Step 1.} Draw a random matrix
$M_k\in \R^{\ell_k\times n}$ 
 from a matrix distribution $\mathcal{M}_k$.
 \\
 \hspace*{1pt} Form a random model $\hathat m_k(\widehat s)$ of the form (\ref{modelGN_S}).
 \\
 \hspace*{1pt}
Compute the inexact step $\widehat s_k$ in (\ref{residuo})$-$(\ref{stepin}).
 Let $s_k=M_k^T \widehat s_k$.
 \item{Step 2.} If $x_k + t_k s_k$ satisfies 
 \be \label{armijo}
f(x_k + t_k s_k) <  f(x_k) + c t_k s_k^T  \nabla f(x_k),
 \ee
 \hspace*{2pt}
 Then (successful iteration)  \\
 \hspace*{15pt} set $  x_{k+1} = x_{k} + t_k s_k$,   $t_{k+1}=\min\{t_{\max},\gamma^{-1} t_k\}$,
$\ell_{k+1}=\max\{\ell_{\min}, \widehat\gamma \ell_k\}$.
 \\
 \hspace*{2pt} Else (unsuccessful iteration)
 \\
 \hspace*{15pt }  set $ x_{k+1} = x_{k}$, $t_{k+1}=\gamma t_k$, $\ell_{k+1}=\min\{\ell_{\max}, \widehat\gamma^{-1} \ell_k\}$.
 \item{Step 3.} Choose $\eta_{k+1}\in [0,  \eta_{\max}]$, $\mu_{k+1}\in [\mu_{\min}, \mu_{\max}]$. Set $ k = k+1$.
 \end{description}
 }
 \label{algo_general}
\vskip 10pt

Successively, in Step 2 of the algorithm, the step search is performed using the Armijo condition (\ref{armijo}) with $c\in (0,1)$ being the small Armijo constant. The test is made on the trial iterate $x_k+t_k s_k$ where $t_k$ is a positive steplength  set at the previous iteration $k-1$. If the test (\ref{armijo}) is satisfied, the iteration is successful, i.e., the trial iterate is accepted, the steplength $t_{k+1}$ is enlarged for the next iteration and the sketching size $\ell_{k+1}$ is reduced for the next iteration taking into account that the current reduced model produced an accepted step.  If the test (\ref{armijo}) fails, the iteration is unsuccessful, i.e., the trial step is discarded, the steplength $t_{k+1}$ is reduced for the next iteration and the sketching size $\ell_{k+1}$ is enlarged. According to a step search strategy, the step direction changes during the backtracking procedure.
Finally, in Step 3 the forcing term $\eta_{k+1}$ and the regularization parameter $\mu_{k+1}$ are defined for the next iteration. 

The use of the Levenberg-Marquardt model instead of the Gauss-Newton model is motivated by the fact that, differently from the Gauss-Newton model, the step  $s_k$ is a descent direction for $f$ as long as the sketched gradient $M_k\nabla f(x_k)=M_k J(x_k)^T F(x_k)$ is nonzero. 
\vskip 5pt
\begin{lemma}\label{discesa}
Let $s_k$ be as in Algorithm \ref{algo_general}. It holds
$
s_k^T\nabla f(x_k)\le -\mu_k\|\widehat s_k\|_2^2, 
$
and  $s_k^T\nabla f(x_k)<0 $ if  $M_k\nabla f(x_k)\neq 0$.
\end{lemma}
\begin{proof}
Let us first consider the case where the system \eqref{nw_sketch} is solved inexactly. Namely, $\widehat s_k$ is computed applying LSMR method and satisfies \eqref{residuo}.
Let us define $G_k=\left(
\begin{array}{c} J(x_k)M_k^T \\ \sqrt{\mu_k} I_{\ell_k}  \end{array} 
\right)$ and $\overline{F}_k= \left(
\begin{array}{c} F(x_k) \\ 0 \end{array} 
\right)$.
Starting from the null initial guess $\widehat s_k^{(0)}=0$, LSMR generates a sequence of iterates
$\{ \widehat s_k^{(j)}\}$, $j\ge 0$, such that
\begin{equation}\label{lsqr}
\textstyle\| G_k \widehat s_k^{(j)}+ \overline{F}_k\|_2^2=\min_{\widehat s\in K_k^{(j)}} \|G_k\widehat s+\overline{F}_k\|_2^2  ,
\end{equation}
with $$K_k^{(j)}=\spa\left\{G_k^T \overline{F}_k, (G_k^T G_k) G_k^T  \overline{F}_k, \ldots, (G_k^T G_k)^{ j-1} G_k^T  \overline{F}_k \right\}.$$
Then, the residual vector $ r_k$ in (\ref{residuo})  corresponding to the inexact step 
 $\widehat s_k=\widehat s_k^{(m)}$, for some $m\ge 0$, is orthogonal  to any vector in $  G_k  K_k^{(m)} $, i.e.,
$\widehat s_k^T G_k^T r_k=\widehat s_k^T \rho_k=0$ with $\rho_k$ defined in \eqref{stepin}. Consequently, \eqref{newin} yields
\begin{eqnarray}
s_k^T \nabla f(x_k) &= & \widehat s_k^T M_k J(x_k)^T F(x_k) \nonumber \\
& =& \widehat s_k^T (-(M_k J(x_k)^T J(x_k)M_k^T+\mu_k I_{\ell_k}) \widehat s_k +\rho_k) \label{inmq}\\
&= &
-s_k^T J(x_k)^T J(x_k) s_k-\mu_k \|\widehat s_k\|_2^2 \le -\mu_k \|\widehat s_k\|_2^2 \nonumber.
\end{eqnarray}
since $J(x_k)^T J(x_k)$ is positive semidefinite. 

In case $\widehat s_k$ solves  the system \eqref{nw_sketch} exactly, the claim follows as above by \eqref{inmq} letting  $\rho_k=0.$

Finally, since by \eqref{newin}-\eqref{stepin}, it holds $\widehat s_k\neq 0$ if and only if $M_k\nabla f(x_k)\neq 0$, the proof is completed.
\end{proof}
\vskip 5pt
Regarding the case where the sketched gradient $M_k \nabla f(x_k)$ is null, we observe that the vector $s_k$ is null and the iteration is unsuccessful. 

The trial step $\widehat s_k\in\R^{\ell_k}$  in (\ref{newin})  gives rise to two relative residuals with respect to the minimization of $\hathat m_k$ and $\widehat m_k$ defined as 
\begin{equation} \label{nuk} 
 \eta_k^*\eqdef\frac{\|\nabla \hathat m_k(\widehat s_k)\|_2}{ \|M_k\nabla f(x_k)\|_2}, \qquad \nu_k^* \eqdef\frac{\|\nabla \widehat m_k(\widehat s_k)\|_2}{\|M_k \nabla f(x_k)\|_2}.
\end{equation}
Such scalars characterize  the approximate solution of the linear systems $\nabla \hathat m_k(\widehat s)=0$ and $\nabla \widehat m_k(s)=0$ respectively.
The next lemma provides a relation between $\eta_k^*$ and $\nu_k^*$ and an upper bound on the norm of $\widehat s_k$
\begin{lemma}\label{boundnuk}
Let us assume that $\widehat s_k\in\R^{\ell_k}$ is such that 
\begin{equation}\label{eq:exactsol}
    (M_kJ(x_k)^T J(x_k)M_k^T + \mu_k I)\widehat s_k = - M_k \nabla f(x_k)+\rho_k.
\end{equation}
and   $\|\rho_k\|_2= \eta_k^* \|M_k\nabla f(x_k)\|_2$. 
Then, if $\eta_k^*=0$
\begin{equation}\label{bound_nustar}
\frac{\mu_k}{\lambda^{1}_k+\mu_k} \le
\nu^*_k\leq\frac{\mu_k}{\lambda^{r_k}_k+\mu_k},
\end{equation}
where  $\lambda^{1}_k$ and $\lambda^{r_k}_k$ are the largest and the smallest nonzero eigenvalue of $M_kJ(x_k)^TJ(x_k)M_k^T$ respectively, otherwise
\begin{equation}\label{gradhatm}
\nu_k^*\leq \frac{\mu_k}{\lambda^{r_k}_k+\mu_k}+ \frac{\lambda_k^1}{\lambda_k^1+\mu_k}\eta_k^*.
\end{equation}
Further, 
\begin{equation}\label{bound_s_hat}
\|\widehat{s}_k\|_2\le \left(\frac{1}{\lambda^{r_k}_k+\mu_k}+\frac{\eta_k}{\mu_k}\right)\|M_k\nabla f(x_k)\|_2.
\end{equation}   
\end{lemma}
\begin{proof}
If $\eta_k^*=0$ then 
$\widehat s_k=-\big(M_kJ(x_k)^T J(x_k)M_k^T + \mu_k I\big)^{-1}M_k \nabla f(x_k)$ and 
\begin{equation}\label{eq:resmodelsketched}
    \nabla \widehat{m}_k(\widehat s_k) =
  \Big(-M_kJ(x_k)^T J(x_k)M_k^T \big(M_kJ(x_k)^T J(x_k)M_k^T + \mu_k I\big)^{-1}+I\Big)M_k \nabla f(x_k) .
\end{equation}
Let $B_k\eqdef M_kJ(x_k)^T J(x_k)M_k^T=Q_k\Lambda_k Q_k^T $ be the eigendecomposition where
$$\Lambda_k = \diag(\lambda^1_k,\dots,\lambda^{r_k}_k,0,\dots,0)\in\R^{\ell_k\times\ell_k},\enspace Q = \big(q^i_k\big|\dots \big|q_k^{\ell_k}\big)\in\R^{\ell_k\times\ell_k},$$
$r_k$ is the rank of the matrix,  $\lambda^1_k\ge \lambda^2_k\ge \cdots,\lambda^{r_k}_k>0$. Note that   $\spa(q_k^1,\dots,q^{r_k}_k) = \range(B_k)$, $\spa(q^{{r_k}+1}_k,\dots q_k^{\ell_k}) = \ker(B_k)$.
Then, 
\begin{equation*}
 \nabla \widehat{m}_k(\widehat s_k) =
Q_k \Big(- \Lambda_k (\Lambda_k+\mu_k I)^{-1} +I\Big)Q_k^T M_k \nabla f(x_k).
\end{equation*}
Since $M_k\nabla f(x_k)\in \range (B_k)$, then $(q_k^i)^T M_k\nabla f(x_k)= 0$ for  $i=r_k+1,\dots \ell_k$, and 
\begin{eqnarray*}
 \nabla \widehat{m}_k(\widehat s_k) &=&  Q_k \Big(-\Lambda_k(\Lambda_k+\mu_k I)^{-1} +I\Big)\left(\begin{matrix}
            (q_k^1)^T M_k\nabla f(x_k)\\
            \vdots\\
    (q_k^{r_k})^T M_k\nabla f(x_k)\\
            0\\
            \vdots\\
            0
        \end{matrix}\right) \\
        &=&- \sum_{i=1}^{r_k}\frac{\mu_k}{\lambda^i_k+\mu_k}\, \Big ( (q_k^i)^T M_k\nabla f(x_k)\Big)q_k^i.
\end{eqnarray*}
The vectors $q_k^i$ are orthonormal,  hence  it follows
\begin{eqnarray*}
\|\nabla \widehat{m}_k(\widehat s_k)\|_2^2 &=&
\left \|\sum_{i=1}^{r_k}\frac{\mu_k}{\lambda^i_k+\mu_k}\, \Big ( (q_k^i)^T M_k\nabla f(x_k)\Big)q_k^i \right \|_2^2\\
&=&
\sum_{i=1}^{r_k} \Big(\frac{\mu_k}{\lambda^i_k+\mu_k}\Big)^2 \, \Big ( (q_k^i)^T M_k\nabla f(x_k)\Big)^2 \|q_k^i\|_2^2 \\
&\le &
 \Big(\frac{\mu_k}{\lambda^{r_k}_k+\mu_k} \Big)^2
\sum_{i=1}^{r_k}  \, \Big ( (q_k^i)^T M_k\nabla f(x_k)\Big)^2  \\
&= &
\Big(\frac{\mu_k}{\lambda^{r_k}_k+\mu_k}\Big)^2
 \| Q^T M_k\nabla f(x_k)\|_2^2  \\
&= &
\Big(\frac{\mu_k}{\lambda^{r_k}_k+\mu_k}\Big)^2
 \| M_k\nabla f(x_k)\|_2^2  \\
\end{eqnarray*}
which gives the upper bound in (\ref{bound_nustar}).  Analogously, the lower bound in (\ref{bound_nustar}) follows.

In the general case $\eta_k^*\ge0$, 
the step $\widehat s_k$ takes the form
\be \label{hatsform}
\widehat s_k = (M_kJ(x_k)^TJ(x_k)M_k^T+\mu_kI)^{-1}(-M_kJ(x_k)^TF(x_k) +\rho_k),
\ee 
and 
taking into account \eqref{stepin}, the equality 
$\|B_k(B_k+\mu_k I)^{-1}\|= \frac{\lambda_k^1}{\lambda^1_k+\mu_k} $ and proceeding as for deriving 
the upper bound, we get (\ref{gradhatm}).

Inequality \eqref{bound_s_hat} is obtained by \eqref{hatsform}
repeating the reasoning above.
\end{proof}

\section{Theoretical analysis}\label{sec3} Algorithm \ref{algo_general} generates a random sequence $\{x_k\}$ since at each iteration the model $\hathat m_k$ is random. Letting $X_k$ be the random variable such that $x_k$ is its realization and $\tau>0$, 
  the hitting time is defined  as
\begin{equation}
 N_{\tau}= \inf \{k: \|\nabla f(X_k)\|_2 \leq \tau \}.\label{hitting_time}
\end{equation}
Following \cite{Shao21}, convergence to a $\tau$-approximate first-order stationary point
 occurs if the algorithm is run for $k\ge  N_{\tau}$ iterations; otherwise the algorithm has not converged.

In this section, exploiting the analysis in 
\cite{cfs} and \cite{Shao21},  we derive a probabilistic bound on the total number of iteration $ N_{\tau}$ required to reach a  $\tau$-approximate first-order stationary point. We perform our analysis making the following assumption on the problem.
\vskip 5pt \noindent
\begin{assumption}\label{regularity}
The objective function $f:\R^n\rightarrow\R$ in problem (\ref{merit}) is continuously differentiable and bounded below by $f_*$. 
The gradient of $f$  is Lipschitz continuous, 
that is, there exist a positive scalar   $L$ such that for any $x,y\in\R^n$
$$\|\nabla f(x)-\nabla f(y)\|_2\leq L\|x-y\|_2.$$
\end{assumption}
\vskip 5pt \noindent
The first requirement for our analysis is to establish some properties of the iterate $x_{k+1}$ when the iteration $k$ is true, namely the random matrix $M_k$  satisfies the following conditions, see \cite{Shao21}.
\vskip 5pt \noindent
\begin{definition}\label{true}
Given  the iteration independent constants $\varepsilon \in (0,1)$, $M_{\max}>0$,
and matrix $M_k\in \R^{\ell_k\times n}$ drawn in Step 	1 of the Algorithm \ref{algo_general},
iteration $k$ is true if 
\begin{eqnarray}
& & \|M_k \nabla f(x_k)\|_2^2\ge (1-\varepsilon)\|\nabla f(x_k)\|_2^2 \label{embedding}\\
& & \|M_k\|_2 \le M_{\max}. \label{bounded_norm}
\end{eqnarray}
\end{definition}
\vskip 5pt
Thus, iteration $k$ is true when $M_k$ provides a one-sided (from below)  embedding with distortion $\varepsilon$ and when its norm is uniformly bounded from above. A relevant property of a true iteration is that it is successful if the steplength $t_k$ is sufficiently small, independently of $k$.
\vskip 5pt \noindent
\begin{lemma}\label{lowertk} Let $x_k$ be the iterate in Algorithm \ref{algo_general} and a true iteration be as in Definition \ref{true}. Suppose that Assumption \ref{regularity}
holds and iteration $k$ is true.
Let $t_{{\rm low}}=\displaystyle \frac{2(1-c)\mu_{\min}}{LM_{\max}^2}$.
If $t_k<t_{{\rm low}}$, then iteration $k$ is successful.
\end{lemma}
\begin{proof}
Using the mean value theorem and Assumptions  \ref{regularity}, we obtain
\begin{eqnarray*}
f(x_{k}+t_{k} s_{k})&=&
f(x_{k}) +  \int_{0}^{1} (\nabla f(x_{k}+w t_{k} s_{k}))^{T} ( t_{k}s_{k}) dw\\
& = & f(x_{k}) +   t_k\nabla f(x_k)^Ts_k+ \int_{0}^{1} t_k(\nabla f(x_{k}+w t_{k} s_{k}) - \nabla f(x_k))^{T} s_{k} dw \\
& \leq & f(x_{k}) + t_k\nabla f(x_k)^Ts_k + \frac{L}{2}t_k^2 \| s_{k}\|_2^2  \\
& \leq & f(x_{k})  + t_k\nabla f(x_k)^Ts_k +  \frac{L}{2}  t_{k}^2 M_{\max}^2\|\widehat s_{k}\|_2^2 .
\end{eqnarray*}
Thus, the Armijo condition (\ref{armijo}) holds  if 
\begin{eqnarray*}
t_k\nabla f(x_k)^Ts_k +  \frac{L M_{\max}^2} {2}  t_{k}^2 \|\widehat s_{k}\|_2^2  < c t_k s_k^T \nabla f(x_k),
\end{eqnarray*}
which is equivalent to $t_k <(c-1) \frac{2}{LM_{\max}^2}   \frac{s_k^T \nabla f(x_k)}{\|\widehat s_k\|^2}$.
Lemma \ref{discesa} implies that $$\frac{2}{L M_{\max}^2} (c-1)  \frac{s_k^T \nabla f(x_k)}{\|\widehat s_k\|^2}> t_{{\rm low}}\eqdef\frac{2 (1-c) \mu_{\min}}{L M_{\max}^2} ,$$ which concludes the proof.
\end{proof}
\vskip 5pt \noindent
Using again the concept of true iteration we can characterize the quantity $f(x_k)-f(x_{k+1})$ for all $k$. We need the following assumption on the generated sequence.
\vskip 5pt 
\begin{assumption}\label{bounded_MJ}
At any true iteration, it holds
$$\sigma_{1}( M_kJ(x_k)^TJ(x_k)M_k^T)\leq \sigma_{\dagger},$$
where $\sigma_{1}(\cdot)$ denotes the maximum singular value  of a  matrix and $\sigma_{\dagger}>0$ is independent of $k$.
\end{assumption}
\vskip 5pt \noindent
\begin{lemma}\label{diff_f_totale} Let $\{x_k\}$ be generated by Algorithm \ref{algo_general}. Let true iterations be defined in Definition \ref{true}.  Suppose that Assumption \ref{bounded_MJ} holds.
\begin{description}
\item{(i)} If iteration $k$ is true and successful with $k<N_{\tau}$, then
$$
f(x_k)-f(x_{k+1})\geq h(\tau, t_k),
$$
where $h:\R^2\rightarrow \R$ is a nonnegative function, non decreasing in its  arguments  $\tau>0$, $t_k>0$.
\item{(ii)}  $f(x_k)- f(x_{k+1})\ge 0$ for all $k\ge 0$. 
\end{description}
\end{lemma}
\begin{proof}
$(i)$ First we prove that
\begin{equation}\label{normgrad}
\|\nabla f(x)\|_2 \le \frac{\sigma_{\dagger}+\mu_{\max}}{(1-\eta_{\max})\sqrt{(1-\varepsilon)}} \|\widehat s_k\|_2.
\end{equation}
In fact, by (\ref{residuo})--(\ref{stepin}),
we have 
\begin{eqnarray*}
\|M_k \nabla f(x_k)\|_2&\le &\|(M_k J(x_k)^T J(x_k)M_k^T+\mu_k I_{\ell_k}) \widehat s_k\|_2+\|\rho_k\|_2 \\
&\le&  (\sigma_{\dagger} +\mu_{\max}) \|\widehat s_k\|_2+\eta_{\max}\|M_k\nabla f(x_k)\|_2.
\end{eqnarray*}
and consequently
$$
\|M_k \nabla f(x_k)\|_2 \le \frac{\sigma_{\dagger} +\mu_{\max}}{1-\eta_{\max}} \|\widehat s_k\|_2.
$$
Hence, using (\ref{embedding}) we obtain (\ref{normgrad}).
Now, if the iteration is true and successful, by  the Armijo condition (\ref{armijo}), Lemma \ref{discesa},
and the inequality (\ref{normgrad})  we have
\begin{eqnarray*}
f(x_{k+1}) &\le & f(x_k) -c t_k\mu_{\min} \|\widehat s_k\|_2^2\\
&\leq&  f(x_k) -c t_k\frac{\mu_{\min}(1-\eta_{\max})^2 (1-\varepsilon)}{( \sigma_{\dagger} +\mu_{\max})^2 }  \|\nabla f(x)\|_2^2 \\
& \le &  f(x_k) -c t_k\frac{\mu_{\min}(1-\eta_{\max})^2 (1-\varepsilon)}{( \sigma_{\dagger} +\mu_{\max})^2 }  \tau ^2.
\end{eqnarray*}
Hence the claim follows with 
\begin{equation}\label{funzh}
h(\tau, t)\eqdef
c t\frac{\mu_{\min}(1-\eta_{\max})^2 (1-\varepsilon)}{(\sigma_{\dagger} +\mu_{\max})^2 }  \tau ^2.
\end{equation}
\vskip 5pt \noindent
$(ii)$ Lemma \ref{discesa}) and the acceptance rule of the step $s_k$ in Algorithm \ref{algo_general} imply   
$f(x_k)- f(x_{k+1})\ge 0$ for all $k\ge 0$. 
\end{proof}
\vskip 5pt
Random matrix distributions guarantee true iterations as in Definition \ref{true} in probability. Following \cite{Shao21}, we suppose that the iterations are true at least with a fixed probability as specified below. 
In what follows, 
$L_k$ is a random variable and $\ell_k$  denotes its realization. 
We make the following assumption.
\vskip 5pt \noindent
\begin{assumption}\label{Atrue}
There exists $\delta_M \in (0,1)$ such that 
\begin{equation}
   \mathbb{P}\left (T_k\,|\,\, X_k=x_k, L_k=\ell_k\right) \geq 1-\delta_M, \quad \quad k=0,1,\ldots,\notag
\end{equation}
where $T_k$ is the event $T_k$=\{iteration $k$ is true\}.
Moreover,  
$T_k$ is conditionally independent on $T_0, T_1, \dots, T_{k-1}$ given $X_k=x_k$  and $L_k=\ell_k$. 
\end{assumption}
\vskip 5pt
Conditions for ensuring the request of true iterations in probability are provided 
in the following Lemma.
\vskip 5pt
\begin{lemma}\label{Lemmatrue}
Suppose that there exist $\varepsilon\in (0,1)$, $\delta_M^{(1)}\in (0,1)$ such that for a(ny) fixed $y\in \{\nabla f(x):\, x\in \R^n\}$,
$M_k\in \R^{\ell_k\times n}$  drawn from a random matrix distribution ${\cal{M}}_k$ satisfies
\be
\mathbb{P} \left(\|M_k  y\|_2^2 \ge (1-\varepsilon)\|y\|_2^2 \, |\, L_k=\ell_k  \right) \ge 1-\delta_M^{(1)}.
\ee
Further, suppose that there exists $\delta_M^{(2)}\in [0,1)$ such that  $M_k$ satisfies
\be
\mathbb{P} \left( \|M_k\|_2 \le M_{\max} \,|\, L_k=\ell_k\right) \ge 1-\delta_M^{(2)},
\ee
where $M_{\max}$ is an iteration independent constant and that  $\delta_M^{(1)}+\delta_M^{(2)}<1$.
\\
For true iterations as in Definition \ref{true}, Assumption \ref{Atrue} is satisfied with $\delta_M=\delta_M^{(1)}+\delta_M^{(2)}$.
\end{lemma}
\noindent
\begin{proof}The proof closely follows that of Lemma 4.4.2 in \cite{Shao21}.    
\end{proof}
\vskip 5pt
In order to derive the result on the hitting time, 
we need the following technical result that represents a generalization to the case of variable sketching size of Lemma A.2 in \cite{cfs}.
 Given $k$, let ${\cal T}_k$ and ${\cal M}_k$ be the random variables corresponding to the realizations $t_k$, $M_k$, respectively and let ${\cal{F}}_{k-1}$ denote the $\sigma$-algebra generated by 
$X_0$, ${\cal T}_0$,  $L_0$, ${\cal M}_0$
$\ldots$, $X_{k-1}$,  ${\cal T}_{k-1}$, $L_{k-1}$, ${\cal M}_{k-1}$, $X_k$, ${\cal T}_k$, $L_k$.
\vskip 5pt
\begin{lemma} \label{lemmaA12cfs}
 Let true iterations be defined in Definition \ref{true}.  
Suppose that Assumption \ref{Atrue} holds with $\delta_M \in (0,1)$.
\begin{description}
\item{i)}
For any $\lambda>0$ and $N\in \mathbb{N}$, we have
\begin{equation}
\E \left[e^{-\lambda \sum_{k=0}^{N-1}  T_k  } \right] \leq 
\left[ e^{( e^{-\lambda}-1)(1-\delta_M)}\right]^N.
\end{equation}
\item{ii)}
If Algorithm \ref{algo_general} runs for $N$ iterations, 
then, for any given $\delta_1 \in (0,1)$,
\begin{equation}\label{}
\mathbb{P}\left( N_S \leq (1-\delta_M)(1-\delta_1)N \right) \leq e^{-\frac{\delta_1^2}{2} (1-\delta_M) N}, 
\end{equation}
where $N_S=\sum_{k=0}^{N-1}  T_k $  is the number of successful iterations.
\end{description}
\end{lemma}
\begin{proof} See the Appendix \ref{Asec3}.
\end{proof}

\vskip 5pt \noindent
\begin{lemma}\label{tmin}
    Given $\tau>0$, suppose that $k<N_{\tau}$ and  that Assumption \ref{bounded_MJ} holds. Let $t_{{\rm low}}$ be defined as in Lemma \ref{lowertk}.
    Then there exist $\psi_{t}= \min\left\{1, \left\lceil\log_{\gamma}\left(\frac{t_{{\rm low}}}{t_0}\right)\right\rceil\right\}\in\mathbb{N}^+$ such that $t_{\min} = t_0\gamma^{\psi_t}$ satisfies $t_{\min}\le \min\{t_{{\rm low}}, \gamma t_0\}$
\end{lemma}
\begin{proof}
    \cite[Lemma 2.1]{cfs}
\end{proof}
\vskip 5pt
We can now state the result on the iteration complexity.
\vskip 5pt
\begin{theorem}\label{th:main}
Suppose that Assumptions \ref{regularity}, \ref{bounded_MJ} and \ref{Atrue} hold  with $\delta_M \in (0,\frac{1}{4})$. 
Let $N_{\tau}$ be defined in \eqref{hitting_time},  $h$ be given in (\ref{funzh}) and $\psi_t$ be given in Lemma \ref{tmin}.
Assume that  Algorithm \ref{algo_general} runs for $N$ iterations. Then, for any
$\delta_1 \in (0,1)$ such that
($1-\delta_M)(1-\delta_1)-\frac 3 4 >0$, if 
\begin{equation}
 N \geq  \left[(1-\delta_M)(1-\delta_1)-\frac 3 4  \right ]^{-1}
\left[\frac{f(x_0)-f_*}{h(\tau, t_0\gamma^{1+\psi_t})} +\frac{\psi_t}{2}\right],
\end{equation}
 we have 
\begin{equation}
  \mathbb{P}\left(N \geq N_{\tau}\right) \geq 1 -e^{-\frac{\delta_1^2}{2}(1-\delta_M)N}. \notag
\end{equation}
\end{theorem}
\begin{proof}
Applying Theorem \cite[Theorem 2.1]{cfs} combined with the results in Lemma \ref{diff_f_totale}, Lemma \ref{lemmaA12cfs} and Lemma \ref{tmin} gives the desired result.  
\end{proof}

\section{Choosing the size $\ell_k$}\label{sec4}
In this section we analyze the step $s_k$ used in the step search strategy, summarize results from the literature on the enforcement of true iterations in probability and introduce a modification to Step 2 of Algorithm \ref{algo_general} that monitors the approximate minimization of the deterministic model $m_k$.

Algorithm \ref{algo_general} can be implemented using random matrix distributions that  generate true iterations in probability according to Definition \ref{true}.
Random ensembles  ${\cal{M}}_k$  which satisfy Lemma \ref{Lemmatrue} are: scaled Gaussian matrices, $s$-hashing matrices,  stable 1-hashing matrices, scaled sampling matrices \cite{acta, Shao21, Wood}. 
For the sake of completeness, in the Appendix \ref{Auno} we report the definition of such distributions and a table summarizing the relations between the  values $\varepsilon$, $M_{\max}$,  $\ell_k$,   $\delta_M^{(1)}$, $\delta_M^{(2)}$.
In principle, due to sketched models, a single iteration of Algorithm \ref{algo_general} is  computationally convenient with respect to the deterministic Levenberg-Marquardt algorithm. But the overall performance of the Algorithm \ref{algo_general} may be worse than that of the deterministic algorithm if the step $s_k$ does not incorporate second-order information from the Gauss-Newton model $m_k$.

Minimizing the reduced model $\widehat m_k$ in (\ref{mhat}) is equivalent to minimizing $m_k$  in the subspace generated by the columns of $M_k^T$. In general, no hint can be given on $s_k = M_k^T\widehat s_k$ as  an approximate minimizer of $m_k$ and on the magnitude of the  scalar  
\begin{equation}\label{thetak}
\theta_k^* \eqdef \frac{\|\nabla  m_k(s_k)\|_2}{\|\nabla f(x_k)\|_2},
\end{equation}
which can be interpreted as a measure of the accuracy of $s_k$ with respect to the optimality condition  $\nabla m_k(s)=0$.
 However, noticing that $\nabla f(x_k)=J(x_k)^T F(x_k)$, this limitation can be overcome  using  a $\varepsilon$-subspace embedding condition for $J(x_k)^T$  and   reformulating the definition of true iteration as follows.

\vskip 5pt
\begin{definition}\label{def:strongtrue}
   Given  the iteration independent constants $\varepsilon \in (0,1)$, $M_{\max}>0$,
and a matrix $M_k\in \R^{\ell_k\times n}$ drawn in Step 	1 of the Algorithm \ref{algo_general},
iteration $k$ is true if 
\begin{eqnarray}
& &  (1+\varepsilon)\|J(x_k)^Tz\|_2^2  \geq  \|M_kJ(x_k)^Tz\|_2^2\geq(1-\varepsilon)\|J(x_k)^Tz\|_2^2 \enspace \text{for every } z\in\R^m,\label{embedding_sub}\\
& & \|M_k\|_2 \le M_{\max}. \label{bounded_norm_sub}
\end{eqnarray}
\end{definition}
\vskip 5pt
Now, recalling the definition of $\eta_k^*$  and $\nu_k^*$ in (\ref{nuk}) we characterize $\theta_k^*$ with respect to $\nu_k$.
\vskip 5pt
\begin{lemma}\label{true_thetastar}
    Let  $\widehat s_k\in\R^{\ell_k}$ be as in (\ref{residuo})--(\ref{stepin}), $s_k=M_k^T\widehat s_k\in \R^n$,  $\nu_k^*$ as in (\ref{nuk}) and   $\theta_k^*$ as in (\ref{thetak}).
Then, if iteration $k$ is true as defined in Definition \ref{def:strongtrue} it holds
$$
\Big(\frac{1-\varepsilon}{1+\varepsilon}\Big) ^{1/2}\nu_k^* \le \theta_k^*\leq \Big(\frac{1+\varepsilon}{1-\varepsilon}\Big) ^{1/2} \nu_k^*.
$$
\end{lemma}
\begin{proof}
Since $ \nabla m_k(s_k)$ and $\nabla f(x_k)$
belong to $\spa (J(x_k)^T)$,  and  $ \nabla \widehat m_k(\widehat s_k)=M_k  \nabla m_k(s_k)$,  the inequality  \ref{embedding_sub} implies
$$
\nu_k^* =   \frac{\|M_k\nabla m_k(s_k)\|_2}{\|M_k \nabla f(x_k)\|_2}\geq \frac{(1-\varepsilon)^{1/2}}{(1+\varepsilon)^{1/2}}\, \frac{\|\nabla m_k(s_k)\|_2}{\|\nabla f(x_k)\|_2} = \Big(\frac{1-\varepsilon}{1+\varepsilon}\Big) ^{1/2}\theta_k^* ,
$$
and
$$
\nu_k^* =   \frac{\|M_k\nabla m_k(s_k)\|_2}{\|M_k \nabla f(x_k)\|_2}\leq \frac{(1+\varepsilon)^{1/2}}{(1-\varepsilon)^{1/2}}\, \frac{\|\nabla m_k(s_k)\|_2}{\|\nabla f(x_k)\|_2} = \Big(\frac{1+\varepsilon}{1-\varepsilon}\Big) ^{1/2}\theta_k^* .
$$
\end{proof}
\vskip 5pt
The property of $\varepsilon$-subspace embedding also yields results on the relation between the rank and the singular values of $M_kJ(x_k)^T$  and  $J(x_k)^T$.
\vskip 5pt
\begin{theorem}Given $\varepsilon\in(0,1)$ and $J(x_k)^T\in\R^{n\times m}$,  suppose that $M_k\in\R^{\ell_k\times n}$ satisfies (\ref{embedding_sub}).
Then 

$i)$ $\rank(M(x_k)J(x_k)^T)=\rank(J(x_k)^T)$ and $\ker(M_kJ(x_k)^T) = \ker(J(x_k)^T)$;

$ii)$ letting $r_k=\rank(J(x_k)^T)$, and
$\sigma_{1}(\cdot)\ge \cdots\ge \sigma_{r_k}(\cdot)$ be the nonsingular values of some given matrix of rank $r_k$, it holds
\begin{equation}\label{bound_sv1}
    \sigma_{1}(M_kJ(x_k)^T)\leq (1+\varepsilon)^{1/2}\sigma_{1}(J(x_k)^T),
\end{equation}
and
\begin{equation}\label{bound_sv}
    \sigma_{r_k}(M_kJ(x_k)^T)\geq (1-\varepsilon)^{1/2}\sigma_{r_k}(J(x_k)^T).
\end{equation}

\end{theorem}
\begin{proof} 
To ease the notation, we drop the iteration index $k$ and we write $M$,  $J$ and $\ell$  in place of $M_k$,  $J(x_k)$ and $\ell_k$ respectively.

The equality $\rank(J^T) =  \rank(MJ^T)$ is proved in \cite[Lemma 2.2.1]{Shao21}. As for the null space of $J^T$ and $MJ^T$, trivially $\ker(J^T)\subseteq \ker(MJ^T)$ holds. Let us assume by contradiction that the inclusion is strict. Then, there exists $\bar z\in\ker(MJ^T)$ such that $\bar z\notin\ker(J^T)$, i.e., 
    $\|MJ^T\bar z\|_2 = 0$ and $\|J^T\bar z\|_2>0$, which contradicts the embedding property \eqref{embedding_sub}.

We now prove the second part of the statement. Let $J^T = U\Sigma V^T$ and $MJ^T = P\widehat{\Sigma} Q^T$ be the singular value decompositions  of $J^T $ and $M J^T$, respectively, with 
$$\Sigma = \left(\begin{matrix} \\   \diag(\sigma_1,\dots,\sigma_r,0,\dots,0)\\
    \\
    O_{(n-m)\times m}\\
    \\
\end{matrix}\right)\in\R^{n\times m},\  \widehat{\Sigma} = \left(\begin{matrix} \\ \diag(\widehat{\sigma}_1,\dots,\widehat{\sigma}_r,0,\dots,0)\\
    \\
    O_{(\ell-m)\times m}\\
    \\
\end{matrix}\right)\in\R^{\ell\times m}$$
where $\sigma_1,\dots,\sigma_r$ and $\widehat{\sigma}_1,\dots,\widehat{\sigma}_r$ are the nonzero singular values of $J^T$ and $MJ^T$, respectively. Moreover
let us denote 
$
V = (v_1,\dots, v_m)\in\R^{m\times m },\ Q = (q_1,\dots, q_m)\in\R^{m\times m } 
$
where $\range(J)=\spa(v_1,\dots,v_r) $, $ \range(JM^T)=\spa(q_1,\dots,q_r)$ and $$\spa(v_{r+1},\dots, v_m) = \ker(J^T) = \ker(MJ^T) = \spa(q_{r+1},\dots,q_m).$$
Note that
\begin{equation}\label{gammar}
    \begin{aligned}
      {\widehat{\sigma}_1^2}&=\|MJ^T q_1\|_2^2\le (1+\varepsilon)\|J^T q_1\|_2^2 = (1+\varepsilon)\|U\Sigma V^T q_1\|_2^2,\\
        \widehat{\sigma}_r^2&= \|MJ^T q_r\|_2^2\geq (1-\varepsilon)\|J^T q_r\|_2^2 = (1-\varepsilon)\|U\Sigma V^T q_r\|_2^2.
    \end{aligned}
\end{equation}
Since $Q$ is orthogonal,  it holds
$$q_r\perp \spa(q_{r+1},\dots,q_m) = \ker(MJ^T) = \ker(J^T)= \spa(v_{r+1},\dots, v_m),
$$
and therefore $q_1^T v_i=q_r^T v_i=0$ for all $i=r+1,\dots,m$.
Thus,
{\begin{equation}
    \begin{aligned}
        \widehat{\sigma}_1^2&  \leq (1+\varepsilon)\sum_{i=1}^r\sigma_i^2(q_1^T v_i)^2\leq (1+\varepsilon)\sigma_1^2\sum_{i=1}^r(q_1^T v_i)^2 \\&= (1+\varepsilon)\sigma_1^2\sum_{i=1}^n(q_1^T v_i)^2 = (1+\varepsilon)\sigma_1^2 \|V^T q_1\|_2^2 = (1+\varepsilon)\sigma_1^2,
    \end{aligned}
\end{equation}}
and
\begin{equation}
    \begin{aligned}
        \widehat{\sigma}_r^2&  \geq (1-\varepsilon)\sum_{i=1}^r\sigma_i^2(q_r^T v_i)^2\geq (1-\varepsilon)\sigma_r^2\sum_{i=1}^r(q_r^T v_i)^2 \\&= (1-\varepsilon)\sigma_r^2\sum_{i=1}^n(q_r^T v_i)^2 = (1-\varepsilon)\sigma_r^2 \|V^T q_r\|_2^2 = (1-\varepsilon)\sigma_r^2,
    \end{aligned}
\end{equation}
which concludes the proof. 
\end{proof}
\vskip 5pt
The previous lemma implies that the subspace embedding property cannot hold if $\ell_k<\rank(J(x_k)^T)$. Further, it characterizes $\lambda^{1}_k$  and $\lambda^{r_k}_k$ in Lemma \ref{boundnuk} since  $\lambda^{1}_k= \sigma_{1}^2(M_kJ(x_k)^T)$  and $\lambda^{r_k}_k= \sigma_{r_k}^2(M_kJ(x_k)^T)$ and thus the condition number $\kappa_2(M_kJ(x_k)^T J(x_k)M_k^T)$ in 2-norm of $M_kJ(x_k)^T J(x_k)M_k^T$ is bounded above by $(1+\varepsilon)/(1-\varepsilon) \kappa_2(J(x_k)^T J(x_k))$.
 
In order to take advantage of random models of reduced dimension, $\ell_k$ should be significantly smaller than $n$.
Scaled Gaussian matrices of dimension $\ell_k\times n$ satisfy (\ref{embedding_sub}) with probability at least $1-\delta$ when $\ell_k={\cal{O}}(\varepsilon^{-2}(r_k+\log(1/\delta)))$ with $r_k$ being the rank of $J(x_k)$, but such matrices are dense and their use is not computationally convenient \cite{Wood}. On the other hand, under suitable conditions, the distribution of $s$-hashing matrices may provide a subspace embedding and computational savings. 
Let us first introduce the notion of coherence $\mu(J(x_k)^T)$ of $J(x_k)^T$.
\vskip 5pt
\begin{definition}  \cite{Mahoney_randalg}
    Given a matrix $J(x_k)^T\in\R^{n\times m}$ with rank $r_k$, let $J(x_k)^T=U_k\Sigma_k V_k^T$ be the economic SVD decomposition where $U_k\in \R^{n\times r_k}$ has  orthonormal columns, $\Sigma_k\in \R^{r_k\times r_k}$ has strictly positive
diagonal entries, $V_k\in \R^{n\times r_k}$ has orthonormal columns. The coherence $\mu(J(x_k)^T)$ of $J(x_k)^T$  is  defined as
    $$\mu(J(x_k)^T)= \max_{i=1, \ldots, n}\|(U_k)_i\|_2,$$
where $(U_k)_i$ denotes the $i$-th row of $U_k$. 
\end{definition}
\vskip 5pt\noindent
It holds $\sqrt{r_k/n}\le \mu(J(x_k)^T)\le 1$, see \cite[Lemma 2.2.3]{Shao21}. 

From the literature we know that  $s$-hashing matrices satisfy \eqref{embedding_sub} in probability under different assumptions on the coherence of the matrix $J(x_k)^T$ and the size $\ell_k.$ In \cite[Theorem 2.3.1]{Shao21} the author proves that  if $\mu(J(x_k)^T)={\cal{O}}(r_k^{-1})$ then 1-hashing matrices satisfy \eqref{embedding_sub} with $\ell_k={\cal{O}}(r_k)$.
 Further,  results  on larger values of  $\ell_k$ state that if $\mu(J(x_k)^T) = O(\log^{-3/2}(r_k))$ then  $\ell_k = O(r_k\log^2(r_k))$ while if no restrictions are put on the coherence, then $\ell_k={\cal{O}}(r_k^2)$ is both necessary and sufficient to enforce a subspace embedding for $1$-hashing matrices (see e.g., \cite[Section 2.3]{Shao21}).  Finally, if $s$-hashing matrices are used with $\ell_k={\cal{O}}(r_k)$, then the request  on $\mu(J(x_k)^T)$ is relaxed by $\sqrt{s}$, see \cite[Theorem 2.4.1]{Shao21}.

We can now draw conclusions on the size of $\ell_k$ based on the discussion above. The overall efficiency of Algorithm \ref{algo_general} depends on two aspects: the use of embedding with small size $\ell_k$ and the rate of convergence since the reduction in the cost of minimizing the model $\hathat m_k$ may be offset by  a large number of iterations performed. At this regard, we observe what follows.
\vskip 5pt
\begin{itemize}
    \item 
Embedding with small size $\ell_k$ can be obtained if the rank of the Jacobian matrix is sufficiently small; this fact occurs if the problem (\ref{merit}) is strongly underdetermined, i.e., $\rank(J(x))\le m \ll n$, or more generally low rank, i.e., $\rank(J(x))\ll \max\{m, n\}$. 
\item 
Though $s_k$ is  a descent direction, see Lemma \ref{discesa}, it may be a poor descent direction for $f$ at $x_k$. 
The rate of convergence   depends on whether $s_k$ retains second-order information from the Gauss-Newton model $m_k$, that ultimately can be measured by the magnitude of $\theta_k^*$.
\item If iteration $k$ is true in the sense of Definition \ref{def:strongtrue}, inequality \eqref{gradhatm} and Lemma \ref{true_thetastar} yield
\be \label{bound_theta}
\theta_k^*\le \Big(\frac{1+\varepsilon}{1-\varepsilon}\Big) ^{1/2} \left(\frac{\mu_k}{\lambda_k^{r_k}+\mu_k}+\frac{\lambda_k^1}{\lambda_k^1+\mu_k}\eta_k^*\right).
\ee
Then, in case of true iterations we have $\theta_k^*=O\left( {\mu_k}/{\lambda_k^{r_k}}+\eta_k^*\right)$
and $s_k$ is an Inexact Gauss-Newton direction corresponding to a forcing term of the order of $\frac{\mu_k}{\lambda_k^{r_k}}+\eta_k^*$. 
\end{itemize}
\vskip 10pt
\algo{}{Revised step 2 of Algorithm \ref{algo_general}. Choosing $\ell_{k+1}$ }{
 \noindent
Given  $c,\, \gamma,\, \widehat\gamma\in (0,1), \,\theta >0$, $t_{\max}$, $\ell_{\min}, \ell_{\max}\in \mathbb{N}$,
$\ell_{\min}<\ell_{\max}\le n$.
\\
Given   $x_k, \, s_k\in {\mathbb R}^n$,  $t_k \in(0, t_{\max}]$, 
$\ell_k \in \mathbb{N}$.
\\
If $x_k + t_k s_k$ satisfies 
$$
f(x_k + t_k s_k) <  f(x_k) + c t_k s_k^T  \nabla f(x_k),
$$
 Then (successful iteration)  \\
\hspace*{23pt} Set $ x_{k+1} = x_{k} + t_k s_k$,   $t_{k+1}=\min\{t_{\max},\gamma^{-1} t_k\}$.
\\
\hspace*{25pt}Compute $\theta_k^*$ in (\ref{thetak}). If 
\begin{equation} \label{modelcond}
   \theta_k^*\le \theta  
\end{equation} 
\hspace*{40pt} 
set $\ell_{k+1}=\max\{\ell_{\min}, \widehat\gamma \ell_k\}$
 \\
\hspace*{25pt}Else 
 \\
\hspace*{40pt} set  $\ell_{k+1}=\min\{\ell_{\max}, \widehat\gamma^{-1} \ell_k\}$.
\\
 \hspace*{2pt} Else (unsuccessful iteration)
 \\
 \hspace*{25pt} set $ x_{k+1} = x_{k}$, $t_{k+1}=\gamma t_k$, $\ell_{k+1}=\min\{\ell_{\max}, \widehat\gamma^{-1} \ell_k\}$.
 }\label{algo_acceptance}
\vskip 10pt

In order to adaptively choose $\ell_k$ and at the same time   to monitor the size of $\theta_k^*$, that in case of false iterations can be large,  we propose a modification in  Step 2 of Algorithm \ref{algo_general} 
as described in Algorithm \ref{algo_acceptance}. We  introduce a prefixed positive threshold $\theta$ and test the magnitude of the value $\theta_k^*$ defined in (\ref{thetak})   with respect to $\theta$. Then, we reduce the sketching size only in case of successful iterations such that $\theta_k^*\le \theta$.
Trivially, setting $\theta=\infty$ inhibits the control (\ref{modelcond}). 
At the extra cost of evaluating the full gradient $\nabla f(x_k)$, Algorithm \ref{algo_acceptance} is a practical procedure for the adaptive choice $\ell_k$ that still allows reductions in the size of the sketching matrices but exploits more information on the current model $m_k$ and step $s_k$ with respect to Algorithm \ref{algo_general}. 
We denote as SLM (Sketched Levenberg-Mar\-quardt) 
the combination of  Algorithms \ref{algo_general}  and \ref{algo_acceptance}.

\section{Local Analysis of the Sketched Levenberg-Marquardt Algorithm}\label{sec5}
In this section we focus on the local convergence of a variant of the SLM Algorithm and show that, despite the use of sketching matrices, in case of true iterations we retain the local error decrease of the deterministic inexact Levenberg-Marquardt approach. 
The procedure is  denoted as SLM-local and  presented  in Algorithm \ref{algo:SLM_loc}. We remark that the steplength $t_k$ is maintained fixed throughout the iterations, i.e., $t_k=1, \, \forall k$. 

\vspace*{10pt}
\noindent
\algo{}{Algorithm SLM-local variant: $k$-th iteration }{
 \noindent
\hspace*{1pt} Given  $\widehat\gamma \in (0,1),\theta >0, \eta_{\max},  \mu_{\max}>0$, $\ell_{\min}, \ell_{\max}\in \mathbb{N}$,
$\ell_{\min}<\ell_{\max}\le n$.
\\
\hspace*{1pt} Given   $x_k\in {\mathbb R}^n$,  $\eta_k\in [0,  \eta_{\max}]$,   $\mu_k\in (0, \mu_{\max}]$,
$\ell_k \in \mathbb{N},\,  \ell_k\in [\ell_{\min},\ell_{\max}]$.
\vskip 2 pt
 \begin{description}
 \item{Step 1.} Draw a random matrix from a matrix distribution $M_k\in \mathcal{M}_k$.
 \\
 \hspace*{1pt} Form a random model $\hathat m_k(\widehat s)$ of the form (\ref{modelGN_S}).
 \\
 \hspace*{1pt}
Compute the inexact step $\widehat s_k$ in (\ref{residuo})$-$(\ref{stepin}).
 Let $s_k=M_k^T \widehat s_k$.
 \item{Step 2.} Set $ x_{k+1} = x_{k} + s_k$.  
\\
\hspace*{2pt} Compute $\theta_k^*$ in (\ref{thetak}). If (\ref{modelcond}) is satisfied then \\
\hspace*{40pt} 
set $\ell_{k+1}=\max\{\ell_{\min}, \widehat\gamma \ell_k\}$
 \\
\hspace*{2pt} Else 
 \\
\hspace*{40pt}  set $\ell_{k+1}=\min\{\ell_{\max}, \widehat\gamma^{-1} \ell_k\}$.
  \item{Step 3.} Choose $\eta_{k+1}\in [0,  \eta_{\max}]$, $\mu_{k+1}\in (0, \mu_{\max}]$. Set $ k = k+1$.
 \end{description}
 }
 \label{algo:SLM_loc}
\vskip 10pt

Let $\Omega^*$ denote the set of all stationary points of $f$, $x^*$ a point in $\Omega^*$, and 
given $\zeta\in(0,1)$, let
$B_\zeta$ be  the closed ball of center $x^*$ and  radius $\zeta$. Moreover, given any $x\in\R^n$, let $\dist(x,\Omega^*)$ denote the distance between $x$ and $\Omega^*,$ i.e. 
$$\dist(x,\Omega^*) = \min\{\|x-z\|_2\ |\ z\in \Omega^*\}.$$

The convergence  analysis  follows the path of  \cite{santos,PILM} where exact and inexact deterministic Levenberg-Marquardt methods are studied and it 
is carried out  under  the following assumptions.

\begin{assumption}\label{ass:lipJ} 
    There exists $L_0>0$ such that for every $x,y\in B_\zeta$ 
    $$\|J(x)-J(y)\|_2\leq L_0\|x-y\|_2.$$
\end{assumption}
\begin{assumption}\label{ass:eigJ} 
For every $x\in B_\zeta$ we have $\rank(J(x)^T J(x))  = \rank(J(x^*)^T J(x^*))= r$, for some positive $r$ and  there exists a positive $\lambda_{\min}$ such that for every $x\in B_\zeta$
$$\min\{\lambda>0 | \lambda\in\eig(J(x)^T J(x))\}\geq \lambda_{\min},$$
where $\eig(J(x)^T J(x))$ denotes the spectrum  of $J(x)^T J(x)$. 
\end{assumption}
\begin{assumption}\label{ass:errorbound}
    There exists $\omega>0$ such that for every $x\in B_\zeta$ 
    $$\omega\dist(x,\Omega^*)\leq \|\nabla f(x)\|_2.$$
\end{assumption}
\begin{assumption}\label{ass:nonlin}
    There exists $\sigma>0$ and $\beta\in[0,1]$ such that for every $x\in B_\zeta$ and every $z\in B_\zeta\cap \Omega^*$ 
    \be \label{bound_ass5.4}
    \|J(x)^T F(z)\|_2\leq \sigma \| x-z\|_2^{1+\beta}.
    \ee
\end{assumption}

Assumption \ref{ass:errorbound} is an error-bound condition with $\|\nabla f(x)\|_2$ as the  residual function \cite{Pang}. This condition is weaker than the full-rank condition, see \cite{santos,noi_errorbound,  PILM, Pang}. 
We note that Assumption \ref{ass:nonlin} is satisfied in case of zero residual problems for any $\sigma\ge 0$. In case of nonzero residual problems it is needed to handle the error due to the employment of $J(x)^TJ(x)$ in place of the true Hessian of $f$.

The following Lemma collects a set of inequalities that follow directly from the Lipschitz continuity of $\nabla f$ and $J$. Note that since $F$ is continuously differentiable, there exists  
${J}_{\max}$ strictly positive such 
$ \|J(x)\|_2\leq {J}_{\max}$ for every $x\in B_\zeta$.
\vskip 5pt
\begin{lemma}\label{lemma_ineq}
Suppose that Assumptions \ref{regularity} and \ref{ass:lipJ} hold. Then
for every $x,y\in B_\zeta$ and every $z\in B_\zeta\cap \Omega^*$,
    \begin{enumerate}
        \item $\|F(y)-F(x)-J(x)(y-x)\|_2\leq L_1\|y-x\|_2^2$, with $L_1 = L_0/2$;
 \item $\|\nabla f(y)-\nabla f(x)-J(x)^T J(x)(y-x)\|_2\leq L_2\|y-x\|_2^2 + \|(J(y)-J(x))^T F(y)\|_2$,\\ with $L_2 = L_1{J}_{\max} $;
        \item $\|(J(y)-J(x))^T F(y)\|_2\leq L_0 {J}_{\max}(\|x-z\|_2\|y-z\|_2+ \|y-z\|_2^2)+ \|J(x)^T F(z)\|_2$\\ \hspace*{110pt} $+\|J(y)^T F(z)\|_2$;
        \item $\|J(x)^T F(x)\|_2\leq L  \dist(x,\Omega^*)$.
    \end{enumerate}
\end{lemma}
\begin{proof}
    See \cite[pp. 1102, 1103]{santos}.
\end{proof}
\vskip 5pt
We now make the following assumption on  the probability of having an iteration true in the sense of Definition \ref{def:strongtrue}. The $\sigma$-algebra $\mathcal{F}_{k-1}$ introduced in Section \ref{sec3} is invoked below.
\vskip 5pt
\begin{assumption}\label{ass:embedding}
Let $\widehat{T}_k$ be the random variable such that
\begin{equation}\label{def:hatT}
\widehat T_k = 
    \begin{cases}
        1 & \text{if \eqref{embedding_sub}-\eqref{bounded_norm_sub} hold at iteration }k    \\
        0 & \text{otherwise.}
    \end{cases}.
\end{equation}
There exists $\pi_M\in(0,1)$ such that for every $k\in\NN$ we have
    $$\P\left(\widehat T_k=1\ |\ \mathcal{F}_{k-1}\right)\geq 1-\pi_M.$$
    Moreover, $\P\left(\widehat T_0=1\right)\geq 1-\pi_M$, and $\widehat T_k$ is conditionally independent on $\widehat T_{k-1},\dots, \widehat T_0$ given $\mathcal{F}_{k-1}.$
\end{assumption}
\vskip 5pt
The following lemma proves that, for true iterations and suitable choices of $\mu_k$ and $\eta_k$, the step $s_k$ is bounded by a multiple of the distance of the current iterate $x_k$ from the set $\Omega^*.$
While in Algorithms \ref{algo_general} and \ref{algo:SLM_loc} the values assigned to the regularization parameter $\mu_k$ and the forcing term $\eta_k$ were  not specified,   here  we  enforce   conditions on the choice of $\mu_k$ and $\eta_k$ in order to recover fast local decrease at true iterations.
\begin{lemma}\label{lemma:loc_s}
Let Assumptions \ref{regularity}, \ref{ass:lipJ}, \ref{ass:eigJ} hold and suppose 
 that there exists a positive constant $\bar c$ such that ${\eta_k}/{\mu_k} =\bar c $ for every $k.$
Then, there exists $c_1$ such that, if $\widehat T_k=1$ and $x_k\in B_\zeta$ then 
$$\|s_k\|_2\leq c_1\dist(x_k,\Omega^*).$$
\end{lemma}
\begin{proof} See the Appendix \ref{Asec5}.
\end{proof}
\vskip 5pt

\begin{lemma}\label{lemma:loc_dist}
    Under  the same assumptions as Lemma \ref{lemma:loc_s} and Assumption \ref{ass:nonlin}, if $x_k,x_{k+1}\in B_\zeta$ and $\widehat T_k=1$ then there exist $L_3,L_4>0$ such that
    $$\omega\dist(x_{k+1},\Omega^*)\leq L_3\dist(x_k,\Omega^*)^{1+\beta}+ L_4 \dist(x_k,\Omega^*)^2 + \theta_k^*L \dist(x_k,\Omega^*),$$
    with $\theta_k^*$ as in (\ref{thetak}).
\end{lemma}
 \begin{proof} See the Appendix \ref{Asec5}.
\end{proof}
\vskip 5pt
\noindent
To proceed in our analysis we let $\bar \eta, \, \bar \mu$ be positive scalars, and
\be \label{bartheta}
    \bar \theta=\frac{(1+\varepsilon)^{1/2}}{(1-\varepsilon)^{1/2}}\left(\frac{\bar \mu}{\lambda_{\min}}+\bar \eta\right),
\ee
with $\varepsilon$ in (\ref{embedding_sub}), and  $\lambda_{\min}$ as in Assumption \ref{ass:eigJ}.
Moreover, for some $\xi\in (0,1)$ let
\begin{equation}
\varsigma =
\left\{
\begin{array}{ll}
\min\left\{ \zeta, \frac{\xi\omega-L_3}{L_4+\bar \theta L^2} \right\} & \mbox { if }  \beta=0,\\
\min\left\{ \zeta, \left(\frac{\xi\omega}{L_3+L_4+\bar \theta L^{1+\beta} }\right)^{1/ \beta} \right \} & \mbox{ if } \beta\in (0,1].
\end{array}
\right. \label{def_varsigma}
\end{equation}
and distinguish the cases  $\beta=0$ and  $\beta \in (0,1]$ in  Assumption \ref{ass:nonlin} and in the following additional assumption.
\vskip 5pt
\begin{assumption}\label{ass: localescalari}
Let $\omega$ be the scalar in Assumption \ref{ass:errorbound}, $\sigma$ be the scalar in Assumption \ref{ass:nonlin}, $c_1$ be the scalar in Lemma \ref{lemma:loc_s}, $L_3,\, L_4$ the scalars in Lemma \ref{lemma:loc_dist}, $\bar \theta$ the scalar in (\ref{bartheta}).
\vskip 2pt
\begin{description}
\item{} If $\beta=0$, suppose that 
$\sigma<\xi\omega/(2+c_1)$, $\eta_k = \bar\eta\|J(x_k)^TF(x_k)\|$, $\mu_k = \bar\mu\|J(x_k)^TF(x_k)\|$.
\vskip 2pt
\item{}
If $\beta \in (0,1]$, suppose that $\eta_k = \bar\eta\|J(x_k)^TF(x_k)\|^\beta$, $\mu_k = \bar\mu\|J(x_k)^TF(x_k)\|^\beta$.
\end{description}
\end{assumption}
\vskip 5pt
We remark that $\eta_k/\mu_k$ is constant as supposed in Lemma \ref{lemma:loc_s}. We also note that in case $\beta=0$, the scalar $\sigma$ is supposed to be sufficiently small. This is in line with the convergence analysis of Gauss-Newton methods for  nonzero residual problems, see \cite{DS}. By definition of $L_3$ (see the proof of lemma \ref{lemma:loc_dist}), the additional condition on  $\sigma$ implies 
$\xi\omega-L_3>0$, and consequently $\varsigma$ in \eqref{def_varsigma} is strictly positive.

We now prove that, assuming that both $x_k$ and $x_{k+1}$ belong to $B_{\zeta}$,  in the true iterations the distance $\dist(x_{k+1},\Omega^*)$ decreases with respect to $\dist(x_{k},\Omega^*)$. Note that, since $\eta_k^*\le \eta_k$ by (\ref{nuk}), then $\theta_k^*$ in (\ref{bound_theta}) satisfies
\be \label{bound_theta2}
\theta_k^*\le \Big(\frac{1+\varepsilon}{1-\varepsilon}\Big) ^{1/2} \left(\frac{\mu_k}{\lambda_k^{r_k}+\mu_k}+\frac{\lambda_k^1}{\lambda_k^1+\mu_k}\eta_k\right).
\ee
\vskip 5pt
\begin{lemma}\label{lemma:loc_xi0}
Let Assumptions \ref{regularity},\ref{ass:lipJ}, \ref{ass:eigJ} and \ref{ass:errorbound} hold. Given any $\xi\in(0,1)$, $\bar\eta\geq 0 $ and $\bar\mu>0$, suppose that  Assumptions \ref{ass:nonlin}  and \ref{ass: localescalari} hold with $\beta\in [0,1]$.
If $\dist(x_{k},\Omega^*)\leq \varsigma$ with $\varsigma$ given in \eqref{def_varsigma},  $x_k,x_{k+1}\in B_{\zeta}$  and $\widehat T_k=1$,
then 
$$\dist(x_{k+1},\Omega^*)\leq\xi\dist(x_k,\Omega^*).$$

\end{lemma}
\begin{proof}
We first consider the case $\beta = 0.$ Inequality \eqref{bound_theta2}, Assumption \ref{ass:eigJ}, the choice 
of $\eta_k$ and $\mu_k$, Item  4 in Lemma \ref{lemma_ineq} yield
\begin{equation}
    \theta_k^*\leq
        \bar\theta L \dist(x_k,\Omega^*),
\end{equation}
with  $\bar \theta$ given in \eqref{bartheta}. 
From Lemma \ref{lemma:loc_dist}, using  $\dist(x_k,\Omega^*)\leq\varsigma$, we have
\be 
  \omega\dist(x_{k+1},\Omega^*)\leq L_3\dist(x_k,\Omega^*)+ L_4 \dist(x_k,\Omega^*)^2 + \bar\theta L^2 \dist(x_k,\Omega^*)^2. 
\ee
This implies
$$
 \omega\dist(x_{k+1},\Omega^*) \le 
 (L_3+ (L_4+\bar\theta L^2 )\varsigma)\dist(x_k,\Omega^*),
 $$
 and  we get the thesis by \eqref{def_varsigma}.
 
In case $\beta\in(0,1]$, inequality \eqref{bound_theta2}, Assumption \ref{ass:eigJ}, the form 
of $\eta_k$ and $\mu_k$, and  Item 4 in Lemma \ref{lemma_ineq} yield
\begin{equation}
    \theta_k^*\leq 
        \bar\theta L^\beta \dist(x_k,\Omega^*)^\beta.
\end{equation}
Using again Lemma \ref{lemma:loc_dist} we get
\be \label{bound_quad}
  \omega\dist(x_{k+1},\Omega^*)\leq L_3\dist(x_k,\Omega^*)^{1+\beta}+ L_4 \dist(x_k,\Omega^*)^2 + \bar\theta L^{1+\beta} \dist(x_k,\Omega^*)^{1+\beta}. 
\ee
Since $\varsigma<1$ by construction, it follows that
$$
\omega\dist(x_{k+1},\Omega^*)\le 
(L_3+ L_4+\bar\theta L^{1+\beta} )\varsigma^\beta\dist(x_k,\Omega^*)
$$
and the thesis follows using \eqref{def_varsigma}.
\end{proof}

\vskip 5pt
We now prove  that if 
$x_{0}$ belongs to $B_{\varsigma}$ and   $\bar k$ consecutive iterations are true, then all the iterates
 $\{x_{k}\}_{k=0}^{\bar k}$ belong to
the ball $B_{\zeta}$ and $\dist(x_k,\Omega^*)$ is smaller than some specified positive scalar.

\vskip 5pt
\begin{lemma}\label{lemma:ind0}
Let Assumptions \ref{regularity},\ref{ass:lipJ}, \ref{ass:eigJ} and \ref{ass:errorbound} hold. Given any $\xi\in(0,1)$, $\bar\eta\geq 0 $ and $\bar\mu>0$, suppose that  Assumptions \ref{ass:nonlin}  and \ref{ass: localescalari} hold with $\beta\in [0,1]$. Let 
\begin{equation}\label{def_barvarsigma}
\bar \varsigma=\min\left\{\varsigma, \frac{\zeta(1-\xi)}{1-\xi+ c_1}\right\}.
\end{equation}
If $x_{0}\in B_{\bar \varsigma}$, and there exists some positive $\bar k$ such that $\widehat T_k=1$ for every $k=0,\dots\bar k$, then we have
    $\dist(x_k,\Omega^*)\leq \bar \varsigma$ and $x_{k+1}\in B_{\zeta}$ for every $k=0,\dots\bar k$.
\end{lemma}
\begin{proof}
    The proof is analogous to that of Lemma 4.2 in \cite{santos}.
\end{proof}
\vskip 5pt
Lemma \ref{lemma:loc_xi0} and \ref{lemma:ind0} above ensure that for all values of $\beta$ the distance of $x_k$ from the set of stationary points decreases at least linearly in case of true iterations.  We now show that the decrease is superlinear and quadratic when $\beta\in(0,1)$ and $\beta=1$, respectively.
\vskip 5pt
\begin{lemma}\label{lemma:locdist2}
    Let the assumptions of Lemma \ref{lemma:ind0}  hold. Let $\beta\in(0,1]$ and $\bar \theta$ given in \eqref{bartheta}. If $x_{0}\in B_{\bar \varsigma}$ and there exists some positive $\bar k$ such that $\widehat T_k=1$ for every $k=0,\dots\bar k$, then we have
    \begin{equation}\label{eq:loc_dist2}
        \dist(x_{k+1},\Omega^*)\leq 
        \frac{L_3+\bar \theta L^{1+\beta} }{\omega}\dist(x_k,\Omega^*)^{1+\beta}+ \frac{L_4}{\omega} \dist(x_k,\Omega^*)^2,
    \end{equation}
    for every $k=0,\dots\bar k$.
\end{lemma}
\begin{proof}
    Since $x_{0}\in B_{\bar \varsigma}$, Lemma \ref{lemma:ind0} ensures that Lemma \ref{lemma:loc_dist} can be applied for all $k=0,\dots,\bar k$.
    Then, \eqref{bound_quad} holds and  gives the thesis.
\end{proof}
\vskip 5pt
The following theorem shows that, in case $\beta=0$ we retain the linear decrease of Levenberg-Marquardt approaches  in $\bar k$ consecutive iterations with probability 
$(1-\pi_M)^{\bar k}$. If Assumption \ref{ass:nonlin} holds with $\beta=1$ we get  quadratic decrease with the same probability.
\vskip 5pt

\begin{theorem}
Let Assumptions \ref{regularity},\ref{ass:lipJ}, \ref{ass:eigJ}, \ref{ass:errorbound} and {\ref{ass:embedding}} hold. Given any $\xi\in(0,1)$, $\bar\eta\geq 0 $ and $\bar\mu>0$, suppose that  Assumptions \ref{ass:nonlin}  and \ref{ass: localescalari} hold with $\beta\in [0,1]$.

If $x_{0}\in B_{\bar \varsigma}$ with $\bar \varsigma$ given in \eqref{def_barvarsigma}, then for every $\bar k\in\mathbb{N}$ we have
$$\P\left(\dist(x_{\bar k},\Omega^*)\leq \xi^{\bar k}\dist(x_{0},\Omega^*)\right )\geq (1-\pi_M)^{\bar k}\quad\quad \mbox{if }\; \beta=0,$$
 and 
$$\P\left(\dist(x_{\bar k},\Omega^*)\leq C^{\sum_{j=0}^{\bar k-1}(1+\beta)^j}\dist(x_{0},\Omega^*)^{(1+\beta)^{\bar k}}\right)\geq (1-\pi_M)^{\bar k}\quad\quad \mbox{if }\; \beta \in (0,1],$$
with $C = (L_3 + L_4 + \bar\theta L^{1+\beta} )/\omega$.
\end{theorem}		
\begin{proof}
First consider the case $\beta=0$ 
and let  $A_j$ be the event $\dist(x_{j+1},\Omega^*)\leq \xi\dist(x_{j},\Omega^*)$.
We prove by induction that 
$$\P\left(\bigcap_{j=0}^{\bar k-1}A_j   \right)\geq(1-\pi_M)^{\bar k},$$
i.e., the probability of linear reduction of $\dist(x_{j+1},\Omega^*)$ with respect to $\dist(x_{j},\Omega^*)$ for $\bar k$ subsequent iteration, $j=0, \ldots,\bar k$, is at least $(1-\pi_M)^{\bar k}$.
For $\bar k=1$, by Lemma \ref{lemma:loc_xi0} we have that $\dist(x_{1},\Omega^*)\leq \xi \dist(x_{0},\Omega^*)$ if $\widehat T_{0}=1$. Therefore,
$$\P(A_0)\geq \P\left(\widehat T_{0} = 1\right)\geq 1-\pi_M.
$$
Let us now assume that 
\begin{equation}\label{prob_dist2}
	\P\left(\bigcap_{j=0}^{\bar k-2}A_j\right)\geq (1-\pi_M)^{\bar k-1}.
\end{equation}
We have
\begin{equation}\label{prob_dist}
    		\begin{aligned}
    			\P\left(\bigcap_{j=0}^{\bar k-1}A_j \right) =\P\left( A_{\bar k-1} \ |\  \bigcap_{j=0}^{\bar k-2}A_j   \right) \P\left(\bigcap_{j=0}^{\bar k-2}A_j\right).
    \end{aligned}
\end{equation}
Let us now consider the first term on the right-hand side of \eqref{prob_dist}. Lemma \ref{lemma:loc_xi0} and \ref{lemma:ind0} ensure that
$$
\P\left ( A_{\bar k-1} \ |\  \bigcap_{j=0}^{\bar k-2}A_j   \right) \geq \P\left ( \widehat T_{\bar k-1} = 1 \ |\  \bigcap_{j=0}^{\bar k-2}A_j   \right)  = \E\left[\widehat T_{\bar k-1} \ |\   \bigcap_{j=0}^{\bar k-2}A_j \right].
$$
Then, by
the definition of expected value and the law of total expectation we have
\begin{equation}\label{prob_dist3}
    		\begin{aligned}
    			\P\left ( A_{\bar k-1} \ |\  \bigcap_{j=0}^{\bar k-2}A_j   \right) 
    			 &\ge \E\left[ \E\left[\widehat T_{\bar k-1} \ |\  \mathcal{F}_{\bar k-2} \right]\ |\   \bigcap_{j=0}^{\bar k-2}A_j \right] \\ &= \E\left[ \P\left(\widehat T_{\bar k-1}=1 \ |\  \mathcal{F}_{\bar k-2} \right)\ |\   \bigcap_{j=0}^{\bar k-2}A_j \right] \\
    			&\geq \E\left[ (1-\pi_M)\ |\   \bigcap_{j=0}^{\bar k-2}A_j \right] = 1-\pi_M.
    		\end{aligned}
\end{equation}
Using inequalities \eqref{prob_dist2} and \eqref{prob_dist3} into \eqref{prob_dist}, we get 
    $$\P\left(\bigcap_{j=0}^{\bar k-1}A_j     \right)\geq(1-\pi_M)^{\bar k}.$$
Therefore
    $$\P\left (\dist(x_{\bar k},\Omega^*)\leq \xi^{\bar k}\dist(x_{0},\Omega^*) \right )\geq \P\left(\bigcap_{j=0}^{\bar k-1}A_j\right)$$
which completes the proof.

If $\beta\in (0,1]$, denoting with $A_j$ the event that inequality \eqref{eq:loc_dist2} holds at iteration $j$,
the proof  follows the above arguments and invokes Lemma \ref{lemma:ind0} and \ref{lemma:locdist2}.
\end{proof}

\vskip 5pt

\section{Numerical results}\label{sec6}
In this section, we investigate the numerical performance of the SLM Algorithm. We also consider the deterministic version  of Algorithm \ref{algo_general} where sketching is not applied, which reduces to the standard Linesearch Levenberg-Marquardt (LLM) procedure  since the direction does not change in case of unsuccessful iterations.

Computing  $s_k$ amounts to using the QR decomposition if $\eta_k=0$, and the LSMR algorithm \cite{LSMR} otherwise. In the case $\eta_k = 0$ we measure the computational cost for solving the linear system as $ 2mn^2+n^2$ and $2m\ell_k^2+\ell_k^2$ for LLM and SLM, respectively (see \cite[Section 2.7.2]{Bjork} for the computational cost of QR with regularization).  If LSMR is applied, the computational cost of the procedure is given by $2m\ell_k q_k$ for the sketched system and $2mn q_k$ for the unsketched system, letting $q_k$ be the number of LSMR iterations performed.

 As for the parameters, we set
 $c = 10^{-4}$,  $\gamma = 0.5$, $\widehat\gamma^{-1}=1.1$, $t_{\max}=1$, 
 $\ell_{\min}=n/10$, $\ell_{\max}=n$,  $\mu_k = 10^{-4}$, $\eta_k=\eta$, $\eta\in [0,1)$, $\forall k\ge 0$.  The maximum number of LSMR iterations is set equal to $\min\{m,\ell_k\}$. The considered matrix distribution ${\cal{M}}_k$ consists of 1-hashing matrices of dimension $\ell_k\times n$. 
 We  denote  SLM$\widehat{p}$ the procedure where the initial sketching size $\ell_0$ is $\widehat{p}\%$  of dimension $n$ and specify the couple $(\eta, \theta)$ used in practice.
We  terminate Algorithms  LLM and SLM  when either $\|\nabla f(x_k)\|_2<10^{-3}$ or $500$ nonlinear iterations are performed. 

\subsection{Problems of varying size and rank}
We  consider a set of artificially generated low rank problems. Given an optimization problem
\begin{equation}\label{lowrankpb1}
   \min_{y\in\R^p} \|\Phi(y)\|_2^2, \quad\text{with}\ \Phi:\R^p\longrightarrow\R^m,
\end{equation}
and a size $n>p$, we 
consider the following augmented problem:
\begin{equation}\label{lowrankpb2}
\min_{x\in\R^n}f(x)=  \|\Phi(Ax)\|_2^2,  
 \end{equation}
where $A$ is random matrix ${p\times n}$ with components uniformly distributed in $[0,1]$ scaled so that $\|A\|_F=1$.
The problems considered are from the CUTEst \cite{CUTEST} collection. 
The number of variables $p$ and observations $m$ of each problem in the collection is determined by a problem-specific parameter $d$. Table \ref{tab:problems1} displays the relation between $m,p$ and $d$ for the considered problems.

\begin{table}
\begin{center}
\begin{tabular}{ c|c|c} 
Problem &  $p$ & $m$\\
\hline
ARTIF & $d+2$ & $d$\\
BRATU2D & $d^2$ & $(d-2)^2$ \\
BROYDN3D & $d$ & $d$ \\
DRCAVTY1& $(d+4)^2$ & $ d^2$ \\
FREURONE& $d$ & $2(d-1)$  \\
OSCIGRNE& $d$ & $d$ \\
\end{tabular}
\end{center}
\caption{CUTEst problems}\label{tab:problems1}
\end{table}

\begin{figure}
\begin{subfigure}[b]{0.47\textwidth}
    \includegraphics[width = 0.9\textwidth]{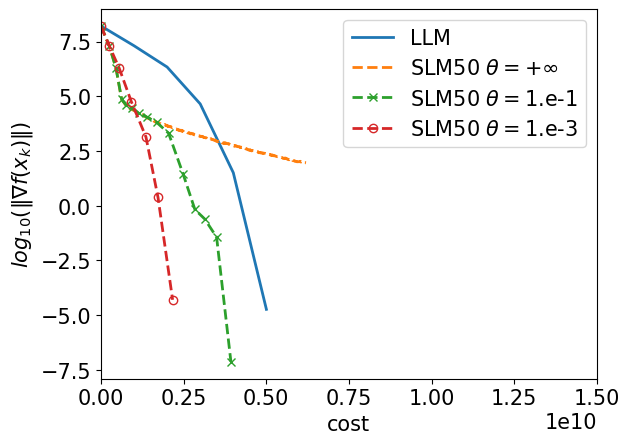}
\end{subfigure}
    \begin{subfigure}[b]{0.47\textwidth}
        \includegraphics[width = 0.9\textwidth]{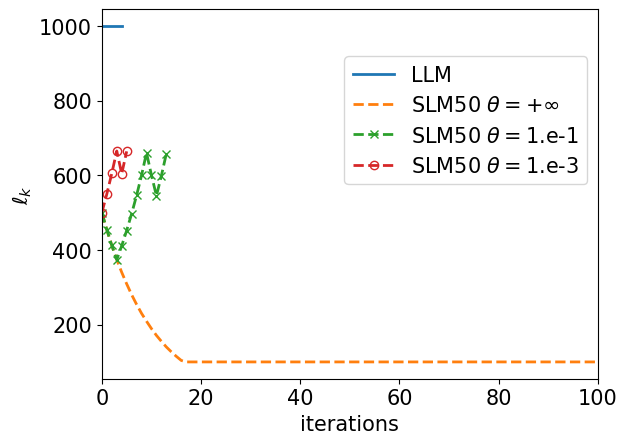}

    \end{subfigure}

    \caption{OSCIGRNE problem $m=500$, $n=1000$. History of LLM and SLM50 $\eta=0$.}
    \label{fig:table_plot}
\end{figure}

To illustrate the effect of using \eqref{modelcond}, we first consider OSCIGRNE problem with $m=500$ and $p=500$, and we define the objective function $f$ as in \eqref{lowrankpb2} with $n=1000$.  We run LLM and SLM50 with $\eta_k=\eta=0,\, \forall k\ge 0$, and initial guess $x_0=(1,\dots,1)^T$.
In Figure \ref{fig:table_plot} we plot
the norm of $\nabla f(x_k)$ versus the computational cost and the value $\ell_k$ versus the iterations. The computational cost is defined as follows. We assign cost $m$ to each evaluation of the residual vector $F(x_k),$ and cost $mn$ to the evaluation of the Jacobian $J(x_k)$. The computational cost for solving the regularized linear system is given by $2m\ell_k^2 + \ell_k^2 $.  In the computation of $\theta_k^*$, the products  $J(x_k)^T F(x_k)$ and    $J(x_k)^T J(x_k)s_k$  cost $3mn$ overall. To summarize, the per-iteration cost is given by $2m\ell_k^2 + \ell_k^2 +4mn + m$.

Figure \ref{fig:table_plot} shows that SLM50 is convenient in terms of cost with respect to LLM when $\theta=10^{-3}$ and $\theta=10^{-1}$.
On the contrary, using $\theta=\infty$ gives poor results. More insight into these results is provided by Table \ref{scalari} and  \ref{scalari_theta}.
For the case where the control \eqref{modelcond} is inhibited, in Table \ref{scalari} we report, for selected iterations $k$,   the values of  $f(x_k)$, $\nabla f(x_k)$, $\ell_k$  along with the relative residuals $\eta_k^*,\nu_k^*,\theta_k^*$, and  the occurrence of successful (S), unsuccessful (U) iteration (Itn).

Table \ref{scalari_theta}  reports the same data obtained with  $\theta = 10^{-1}$ in \eqref{modelcond}.

Table \ref{scalari}  shows that 
the reported iterations are successful and consequently $\ell_k$ decreases steadily, as shown in Figure \ref{fig:table_plot} (right). But the  decrease of $f(x_k)$ and $\|\nabla f(x_k)\|_2$ is low starting from $k=4$, and this behavior can be ascribed to the the quality of $s_k$ with respect to the model $m_k$ in (\ref{mk}). The values of $\eta_k^*$ are nonzero but very small due to round-off precision, the values $\nu_k^*$ are small as expected, while the value of $\theta_k^*$ are close to one for $k\ge 4$. This occurrence affects the convergence history and  SLM50 compares poorly with LLM. On the contrary, Table \ref{scalari_theta} shows that imposing the control (\ref{modelcond}) gives rise to an increase in the size $\ell_k$ at some iterations, even if the iteration is successful, and  greatly improves the performance of SLM50.

\begin{table}
\begin{center}
\begin{tabular}{ c|c|c|c|c|c|c|c} 
$k$ &  $f(x_k)$ & $\nabla f(x_k)$ & $\ell_k$& $\eta_k^*$ & $\nu_k^*$ &$\theta_k^*$  & Itn\\
\hline
0 & 3.50e$+$8& 1.65e$+$8 &500& 2.41e$-$16 & 1.83e$-$11 & 1.93e$-3$& S\\
1 & 1.73e$+$7& 2.06e$+$7 &454& 2.61e$-$16 & 8.18e$-$11 & 4.54e$-$3& S\\
2 & 4.73e$+$5& 1.94e$+$6 &412& 5.13e$-$16 & 4.17e$-$10 & 2.65e$-$2& S\\

3 & 1.10e$+$5& 7.37e$+$4 &374& 1.03e$-$14 & 5.72e$-$9 & 5.46e$-$1& S\\
4 & 6.97e$+$4& 4.01e$+$4 &340& 8.10e$-$15 & 8.17e$-$9 & 7.42e$-$1& S\\
5 & 4.70e$+$4& 2.98e$+$4 &309& 6.56e$-$15 & 5.52e$-$9 & 8.82e$-$1& S\\

6 & 3.69e$+$4& 2.63e$+$4 &280& 6.61e$-$15 & 5.34e$-$9 & 9.02e$-$1& S\\
7 & 2.99e$+$4& 2.37e$+$4 &254& 4.93e$-$15 & 4.37e$-$9 & 9.53e$-$1& S\\
8 & 2.38e$+$4& 2.26e$+$4 &230& 3.96e$-$15 & 2.60e$-$9 & 8.29e$-$1& S\\

9 & 2.02e$+$4& 1.87e$+$4 &209& 2.86e$-$15 & 1.88e$-$9 & 9.16e$-$1& S\\
10 & 1.84e$+$4& 1.87e$+$4 &189& 3.54e$-$15 & 1.83e$-$9 & 9.11e$-$1& S\\
100 & 1.82e$+$3& 3.93e$+$3 &100& 1.48e$-$15 & 6.03e$-$10 & 9.74e$-$1& S\\
200 & 2.40e$+$2& 1.36e$+$3 &100& 1.30e$-$15 & 5.23e$-$10 & 9.79e$-$1& S\\
300 & 3.83e$+$1& 5.10e$+$2 &100& 1.10e$-$15 & 5.32e$-$10 & 1.01e$+$0 & S\\
400 & 7.35e$+0$& 2.30e$+$2 &100& 1.27e$-$15 & 6.44e$-$10 & 1.02e$+$0& S
\end{tabular}
\end{center}
\caption{OSCIGRNE problem   $m=500$, $n=1000$. History of  SLM50 along the iterations,  $\eta= 0$, $\theta = +\infty$. Itn: successful iteration (S), unsuccessful iteration (U).
}\label{scalari} 
\end{table}

\begin{table}
\begin{center}
\begin{tabular}{ c|c|c|c|c|c|c|c} 
$k$ &  $f(x_k)$ & $\nabla f(x_k)$ & $\ell_k$& $\eta_k^*$ & $\nu_k^*$ &$\theta_k^*$  & Itn\\
\hline
0 & 3.50e$+$8& 1.64e$+$8 &500& 4.20e$-$16 & 2.35e$-$11 & 1.54e$-$3& S\\
1 & 1.72e$+$7& 2.07e$+$7 &454& 2.43e$-$16 & 6.57e$-$11 & 4.38e$-$3& S\\
2 & 4.53e$+$5& 1.94e$+$6 &412& 5.90e$-$16 & 4.19e$-$10 & 2.46e$-$2& S\\

3 & 9.44e$+$4& 7.05e$+$4 &374& 6.07e$-$15 & 5.26e$-$9 & 4.98e$-$1& S\\
4 & 6.06e$+$4& 3.50e$+$4 &411& 1.53e$-$14 & 1.09e$-$8 & 8.43e$-$1& S\\
5 & 3.81e$+$4& 2.96e$+$4 &452& 3.88e$-$14 & 2.11e$-$8 & 5.91e$-$1& S\\

6 & 1.76e$+$4& 1.75e$+$4 &497& 2.82e$-$14 & 2.89e$-$8 & 7.11e$-$1& S\\
7 & 8.04e$+$3& 1.25e$+$4 &546& 4.35e$-$14 & 3.33e$-$8 & 5.23e$-1$& S\\
8 & 2.08e$+$3& 6.54e$+$3 &600& 2.09e$-$13 & 1.48e$-$7 & 3.07e$-$1& S\\

9 & 2.24e$+$2& 2.02e$+$3 &660& 9.74e$-$14 & 7.76e$-$8 & 1.42e$-$6& S\\
10 & 1.64e$-$4& 2.84e$+$1 &600& 6.05e$-$15 & 5.69e$-$9 & 9.43e$-$3& S\\
11 & 3.51e$-$6& 2.67e$-$1 &545& 2.82e$-$14 & 3.05e$-$8 & 5.08e$-$1& S\\

12 & 1.14e$-$6& 1.36e$-$1 &599& 8.75e$-$14 & 8.06e$-$8 & 3.73e$-$1& S\\
13 & 1.29e$-$7& 5.06e$-$2 &658& 9.08e$-$14 & 8.26e$-$8 & 1.67e$-$6 & S\\
14 & 5.0e$-$19& 8.67e$-$8 & & & & &  \\

\end{tabular}
\end{center}
\caption{OSCIGRNE problem  $m=500$, $n=1000$. History of LLM and SLM50 along the iterations,   $\eta=0$, $\theta = 10^{-1}$. Itn: successful iteration (S), unsuccessful iteration (U).}\label{scalari_theta} 
\end{table}

\begin{figure}[h]
\centering
    \begin{subfigure}[b]{0.327\textwidth}
    \centering
        \includegraphics[width = \textwidth]{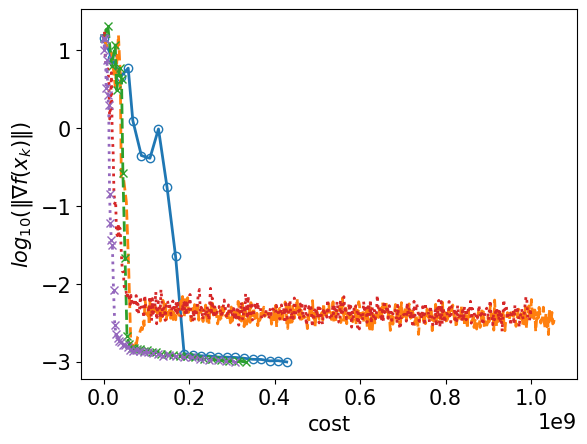}
 \caption{ARTIF}
    \end{subfigure}
    \begin{subfigure}{0.327\textwidth}
    \centering
        \includegraphics[width = \textwidth]{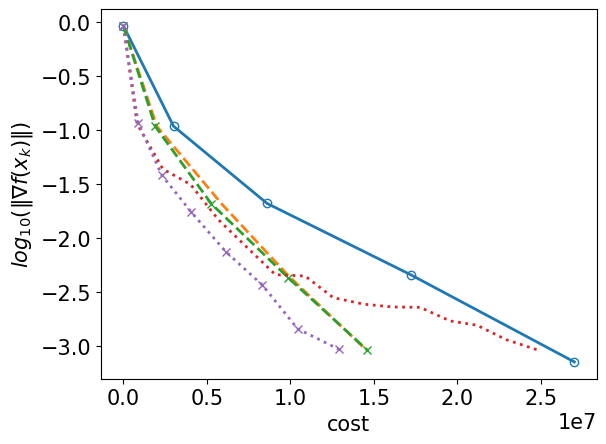}
 \caption{BRATU2D}
    \end{subfigure}
       \centering
        \begin{subfigure}{0.327\textwidth}
      \centering
        \includegraphics[width = \textwidth]{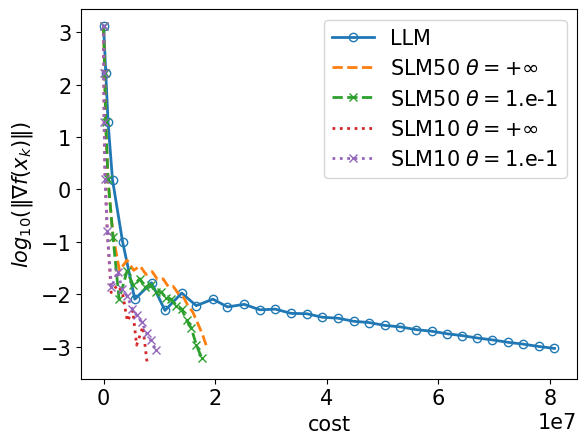}
 \caption{BROYDN3D}
    \end{subfigure}

    \begin{subfigure}[b]{0.327\textwidth}
    \centering
        \includegraphics[width = \textwidth]{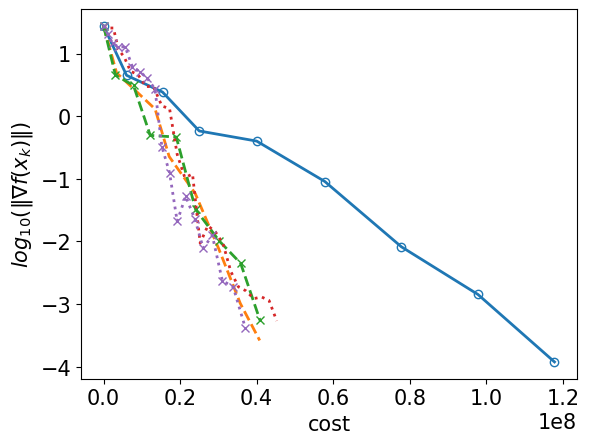}
 \caption{DRCAVTY1}
    \end{subfigure}
    \begin{subfigure}{0.327\textwidth}
    \centering
        \includegraphics[width = \textwidth]{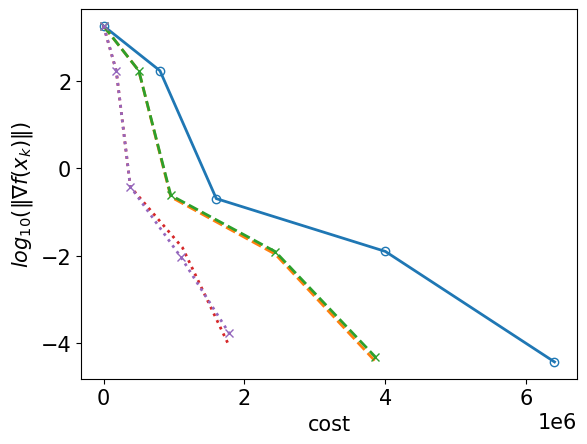}
 \caption{FREURONE}
    \end{subfigure}
       \centering
        \begin{subfigure}{0.327\textwidth}
      \centering
        \includegraphics[width = \textwidth]{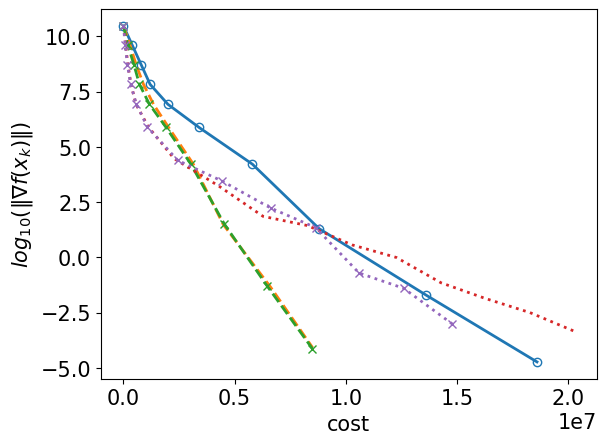}
 \caption{OSCIGRNE}
    \end{subfigure}
   
    \caption{ CUTEst problems solved with LLM, SLM10 and SLM50, $\eta = 10^{-3}, \theta \in\{ +\infty, 10^{-1}\}$. Norm of $\nabla f(x_k)$ vs computational cost.}\label{fig:cute100_grad}

    \end{figure}

\begin{figure}[h]
\centering
    \begin{subfigure}[b]{0.327\textwidth}
    \centering
        \includegraphics[width = \textwidth]{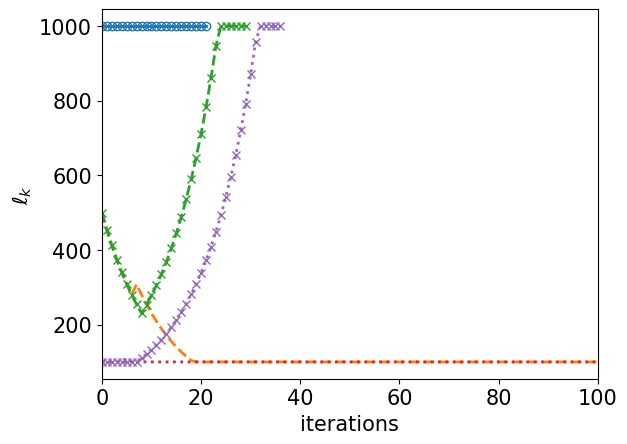}
 \caption{ARTIF}
    \end{subfigure}
    \begin{subfigure}{0.327\textwidth}
    \centering
        \includegraphics[width = \textwidth]{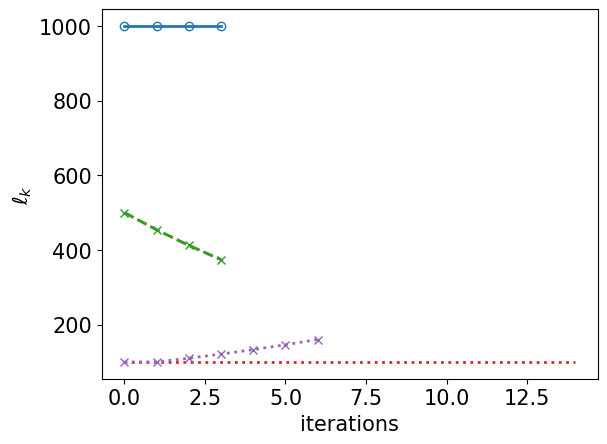}
 \caption{BRATU2D}
    \end{subfigure}
       \centering
        \begin{subfigure}{0.327\textwidth}
      \centering
        \includegraphics[width = \textwidth]{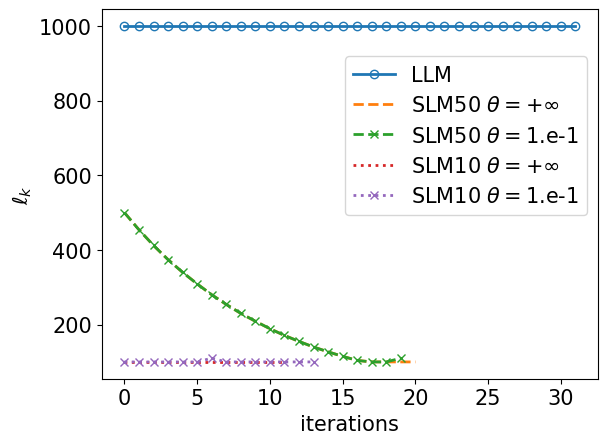}
 \caption{BROYDN3D}
    \end{subfigure}

    \begin{subfigure}[b]{0.327\textwidth}
    \centering
        \includegraphics[width = \textwidth]{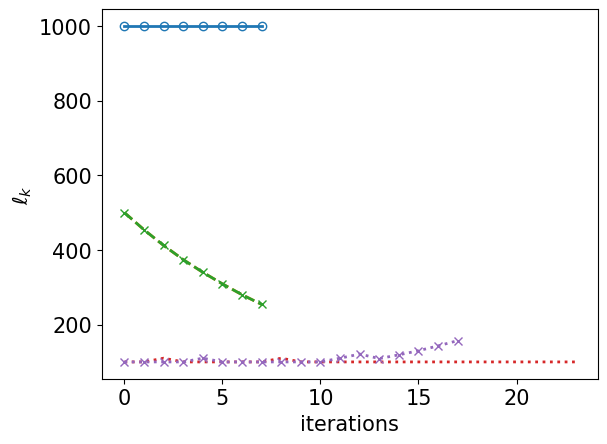}
 \caption{DRCAVTY1}
    \end{subfigure}
    \begin{subfigure}{0.327\textwidth}
    \centering
        \includegraphics[width = \textwidth]{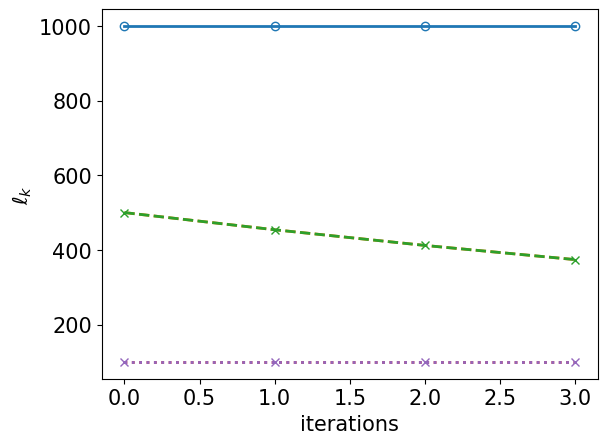}
 \caption{FREURONE}
    \end{subfigure}
       \centering
        \begin{subfigure}{0.327\textwidth}
      \centering
        \includegraphics[width = \textwidth]{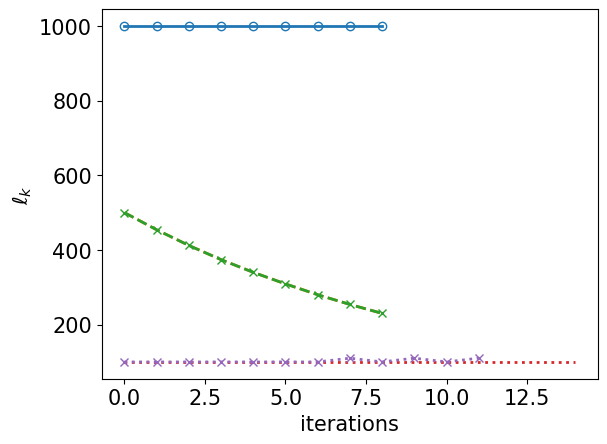}
 \caption{OSCIGRNE}
    \end{subfigure}
   
    \caption{ CUTEst problems solved with LLM, SLM10 and SLM50, $\eta = 10^{-3}, \theta \in\{ +\infty, 10^{-1}\}$.  Sketching size $\ell_k$  vs iterations.}\label{fig:cute100_size}

    \end{figure}

Further experiments were carried out solving strongly underdetermined problems where $m=100$ and the artificial size $n$ is set to $1000$. 
The initial guess was fixed as $x_0=(1,\dots,1)^T$ and  SLM$p$ algorithm was run 11 times for each tested value $p$. The LLM method is deterministic, and therefore it is not necessary to run it repeatedly.
We present results  obtained with the constant forcing terms $\eta=10^{-3}$, $\theta\in \{\infty,  10^{-1}\}$,  and plot the median one in terms of overall computational cost.  The computational cost is measured as described above taking into account that  each iteration of LSMR method has cost $2m\ell_k$. Hence, the per-iteration cost is given by
$2m\ell_k q_k+4mn + m$,
where $q_k$ is the number of LSMR iterations performed at iteration $k$. 
In Figure \ref{fig:cute100_grad} the norm of the gradient is plotted against the cost, in Figure \ref{fig:cute100_size} the subspace dimension $\ell_k$ is displayed versus the iterations.  

The sketched algorithms with $\theta=10^{-1}$ perform well compared to LLM algorithm  on all the considered problems. In the solution of  DRCAVTY1 problem, SLM10 and SLM50 perform significantly better than LLM  for all $\theta.$ In the solution of BROYDN3D and FREURONE problems, SLM10 algorithm  appears to be the most effective and the value $\theta$ used does not  affect the performance significantly. For problem OSCIGRNE SLM50 is significantly cheaper than LLM and is not affected by the choice of $\theta,$ while SLM10 is comparable to LLM for $\theta = 10^{-1}$ and more costly for $\theta=+\infty.$ In the solution of problem BRATU2D,  SLM50 with both choices of $\theta$ and SLM10 with $\theta = 10^{-1}$  perform similarly and significantly better than LLM, while the method that employs $\ell_0=10$ and $\theta = +\infty$ is comparable to LLM. ARTIF is the only considered problem where employing the sketching does not seem to yield significant gains in terms of computational cost. If \eqref{modelcond} is inhibited, runs fails  while  using $\theta=10^{-1}$ yields to $\ell_k=n$ in the last iterations, as shown in Figure \ref{fig:cute100_size}. For the other problems, Figure \ref{fig:cute100_size} shows that the sketching size $\ell_k$ is truly adaptive and changes along the iterations. It's behavior heavily depends on the  parameters $\theta, \ell_0$ employed and on the specific  problem.  We note that some curves of SLM$\widehat{p}$ are overlapped for a given $\widehat{p}$ and varying $\theta$.

\begin{figure}[h]
\centering
    \begin{subfigure}[b]{0.327\textwidth}
    \centering
        \includegraphics[width = \textwidth]{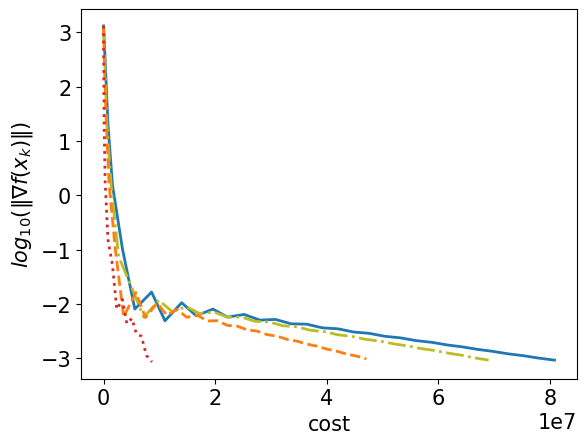}
 \caption{BROYDN3D}
    \end{subfigure}
    \begin{subfigure}{0.327\textwidth}
    \centering
        \includegraphics[width = \textwidth]{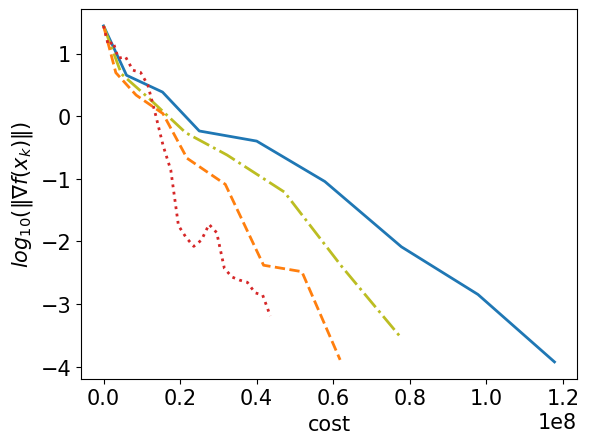}
 \caption{DRCAVTY1}
    \end{subfigure}
       \centering
        \begin{subfigure}{0.327\textwidth}
      \centering
        \includegraphics[width = \textwidth]{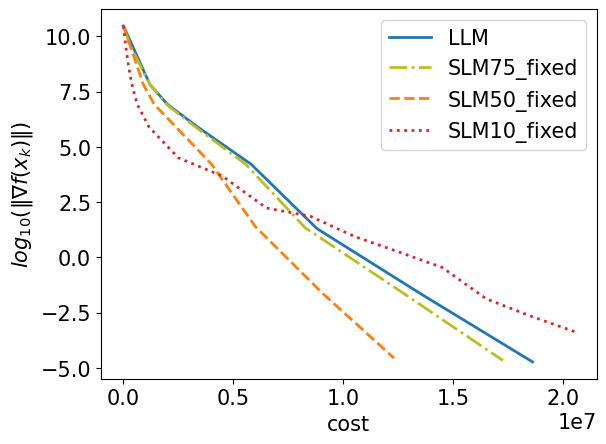}
 \caption{OSCIGRNE}
    \end{subfigure}

    \caption{ CUTEst problems solved with $\eta = 10^{-3}$ and constant $\ell_k = \ell_0$, $\forall k$. Results for LLM, SLM75\_fixed, SLM50\_fixed and SLM10\_fixed. Norm of $\nabla f(x_k)$ vs computational cost.}\label{fig:cute100_grad_fixed}

    \end{figure}

We conclude the current set of experiments showing that  the adaptive choice of the sketching size positively affects the performance of the algorithm SLM. Hence, we repeat the tests conducted above using constant sketching dimension,  i.e.,  $\ell_k=\ell_0$, $\forall k$. All other parameters are set as in the previous tests. We solved BROYDN3D, DRCAVTY1 and OSCIGRNE, $m=100$, $n=1000$,  with  $\ell_0\in\{750,500,100\}$, corresponding to $75\%, 50\%$ and $10\%$ of the problem  dimension, respectively. We denote SLM$\widehat{p}$\_fixed the corresponding algorithm.  The results are reported in Figure \ref{fig:cute100_grad_fixed} and we are interested in comparing such results with those in Figure \ref{fig:cute100_grad} obtained with the adaptive strategy. We already noticed in Figure \ref{fig:cute100_size} that, on the considered problems, for SLM10 with $\theta = +\infty$, the size $\ell_k$ remains constant and always equal to 100; therefore  the behavior of SLM10 and SLM10\_fixed are analogous. 
Regarding the fixed values $\ell_k = 500$ and $\ell_k = 750$, $\forall k$, the performance of SLM75\_fixed is significantly worse than that of SLM50\_fixed and SLM50, and  in two problems out of three, the overall cost of SLM75\_fixed is comparable to that of LLM. 
Moreover, for the three considered problems, the final computational cost of the SLM50\_fixed is significantly higher than the cost  of SLM50.
Overall the results of Figures \ref{fig:cute100_grad}-\ref{fig:cute100_grad_fixed} suggest that SLM with constant sketching size can work well in practice, but the most effective sketching size seems to depend heavily on the considered problem. Hence, employing an adaptive strategy for the choosing $\ell_k$ seems to improve the robustness and the performance of the SLM strategy.

\begin{figure}[h]
\centering
    \begin{subfigure}[b]{0.327\textwidth}
    \centering
        \includegraphics[width = \textwidth]{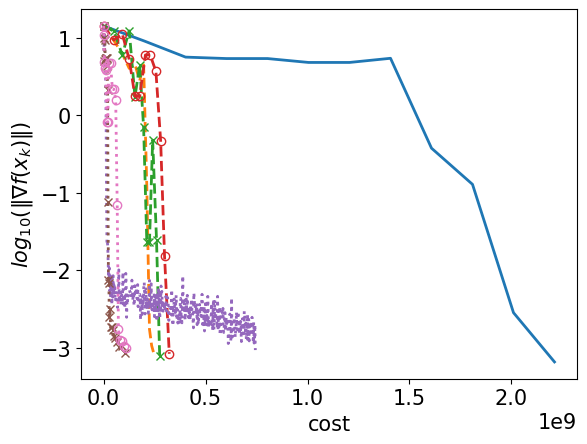}
 \caption{ARTIF}
    \end{subfigure}
    \begin{subfigure}{0.327\textwidth}
    \centering
        \includegraphics[width = \textwidth]{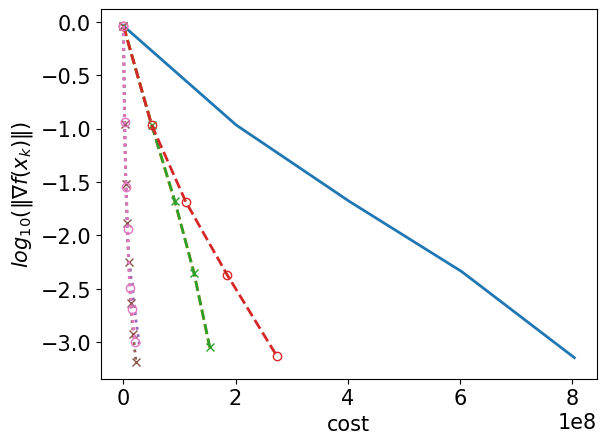}
 \caption{BRATU2D}
    \end{subfigure}
       \centering
        \begin{subfigure}{0.327\textwidth}
      \centering
        \includegraphics[width = \textwidth]{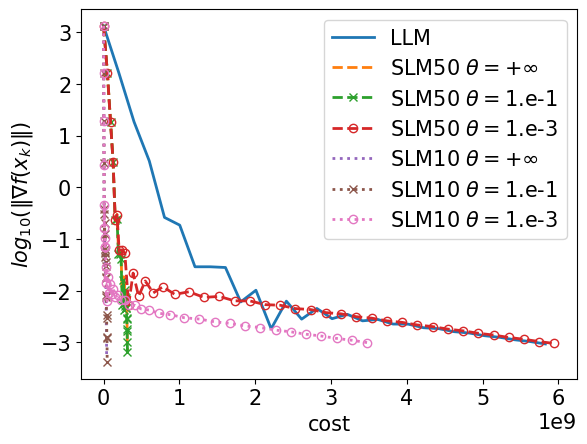}
 \caption{BROYDN3D}
    \end{subfigure}

    \begin{subfigure}[b]{0.327\textwidth}
    \centering
        \includegraphics[width = \textwidth]{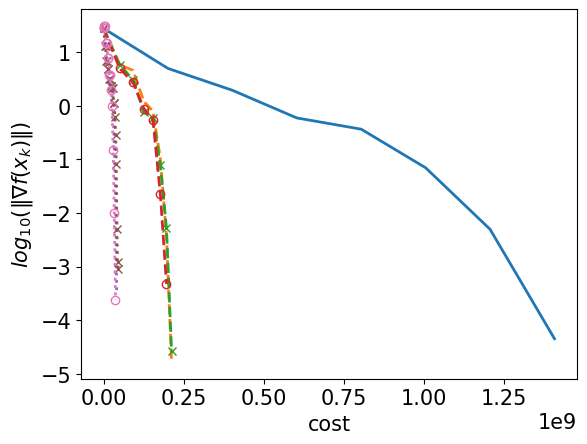}
 \caption{DRCAVTY1}
    \end{subfigure}
    \begin{subfigure}{0.327\textwidth}
    \centering
        \includegraphics[width = \textwidth]{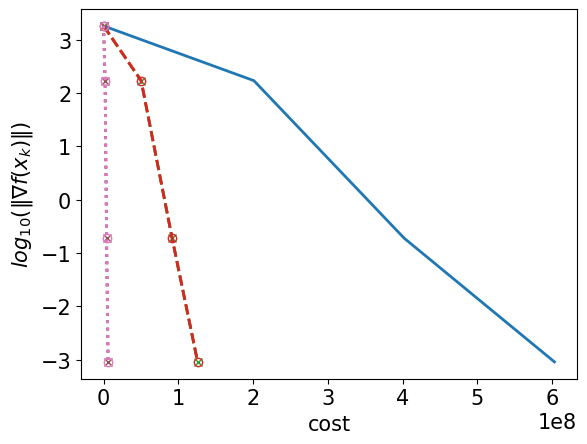}
 \caption{FREURONE}
    \end{subfigure}
       \centering
        \begin{subfigure}{0.327\textwidth}
      \centering
        \includegraphics[width = \textwidth]{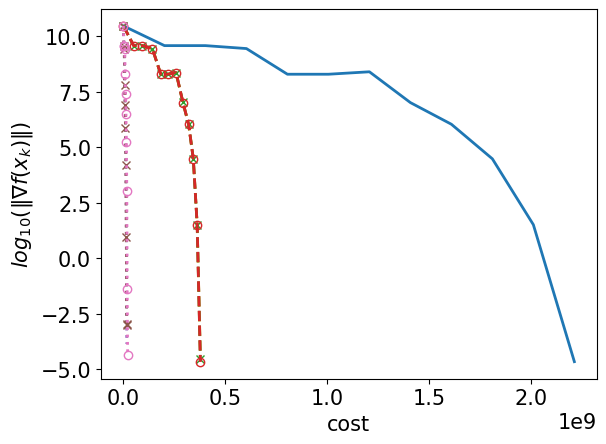}
 \caption{OSCIGRNE}
    \end{subfigure}

    \caption{ CUTEst problems solved with LLM, SLM10 and SLM50, $\eta=0$, $ \theta \in \{ +\infty, 10^{-1}, 10^{-3}\}$.  Norm of $\nabla f(x_k)$ vs computational cost.}\label{fig:cute100_exact}

    \end{figure}

\begin{figure}[h]
\centering
    \begin{subfigure}[b]{0.327\textwidth}
    \centering
        \includegraphics[width = \textwidth]{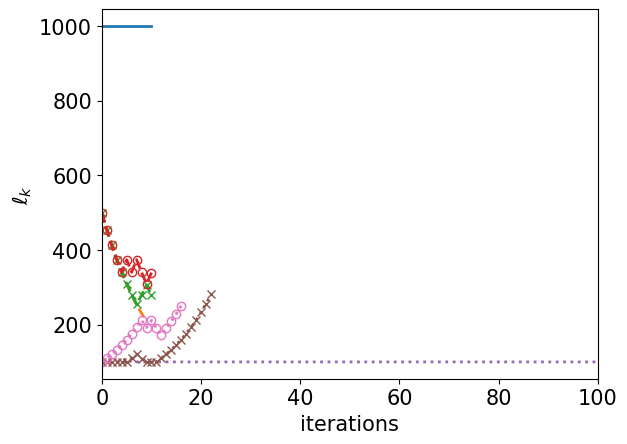}
 \caption{ARTIF}
    \end{subfigure}
    \begin{subfigure}{0.327\textwidth}
    \centering
        \includegraphics[width = \textwidth]{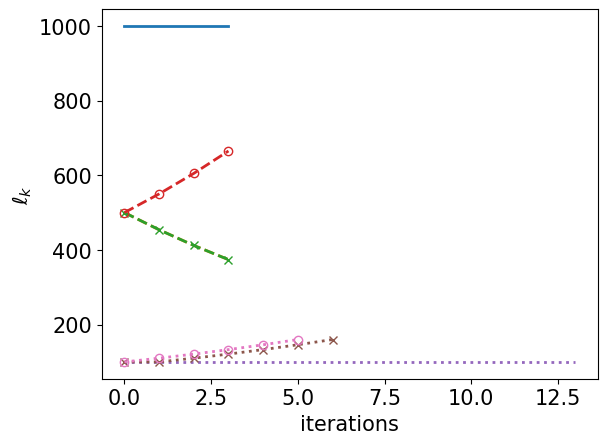}
 \caption{BRATU2D}
    \end{subfigure}
       \centering
        \begin{subfigure}{0.327\textwidth}
      \centering
        \includegraphics[width = \textwidth]{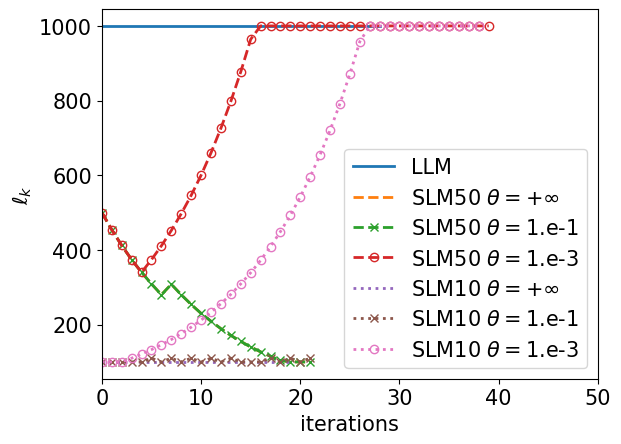}
 \caption{BROYDN3D}
    \end{subfigure}

    \begin{subfigure}[b]{0.327\textwidth}
    \centering
        \includegraphics[width = \textwidth]{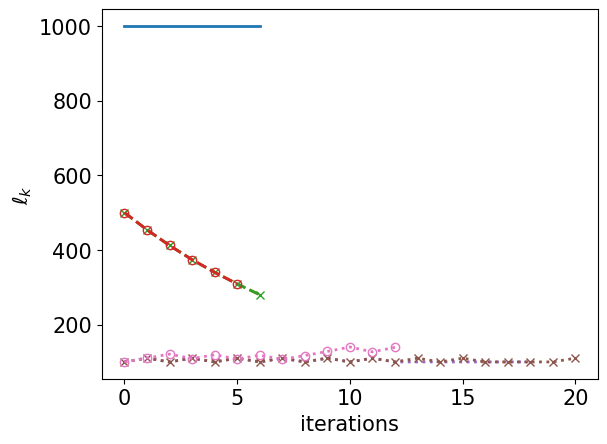}
 \caption{DRCAVTY1}
    \end{subfigure}
    \begin{subfigure}{0.327\textwidth}
    \centering
        \includegraphics[width = \textwidth]{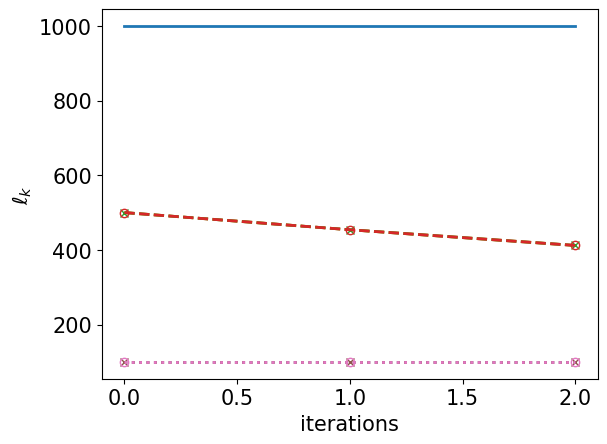}
 \caption{FREURONE}
    \end{subfigure}
       \centering
        \begin{subfigure}{0.327\textwidth}
      \centering
        \includegraphics[width = \textwidth]{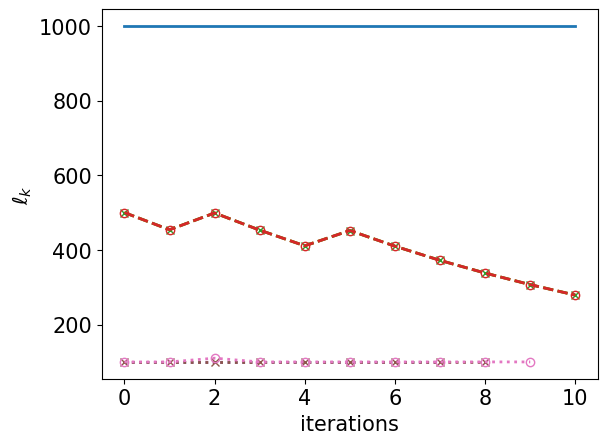}
 \caption{OSCIGRNE}
    \end{subfigure}

    \caption{ CUTEst problems solved with LLM, SLM10 and SLM50, $\eta=0$, $ \theta \in \{ +\infty, 10^{-1}, 10^{-3}\}$. Embedding size $\ell_k$ vs iterations.}\label{fig:cute100_exact_size}

    \end{figure}

We repeated the tests in Figures \ref{fig:cute100_grad}, \ref{fig:cute100_size} using a direct method for the solution of the linear system, $\eta=0$. 
The results are reported in Figures \ref{fig:cute100_exact} and \ref{fig:cute100_exact_size}. Figure \ref{fig:cute100_exact} shows that SLM10 and SLM50 procedures are effective in all runs except the case when $\theta=10^{-3}$ is used  for problem BROYDN3D. In such costly run,   Figure \ref{fig:cute100_exact_size} displays that the condition (\ref{modelcond}) is not satisfied and the embedding size increases steadily.

Summarizing the results presented, SLM$\widehat{p}$ showed to be effective in terms of computational cost. For some problems, testing the condition (\ref{modelcond}) was crucial for improving the performance of SLM$p$ algorithm; further using moderate values of $\theta$, such as $\theta=10^{-1}$ did not deteriorate the behavior of the sketched algorithm.

\subsection{Binary Classification}
We consider a binary classification problem  with logistic model and least-squares loss of the form (\ref{merit}) where
\begin{equation}\label{classification}
    F_{i}(x) = b_i-\frac{1}{1+e^{-x^T a_i}}, \quad i=1,\dots,m,\end{equation}
and  $a_i\in \R^n$, $b_i\in\{0,1\}$  are the features vectors and the labels of the training set   respectively. 

The datasets used are  GISETTE \cite{gisette}  and  REJAFADA  \cite{rejafada}. 
Regarding GISETTE, the problem dimension is $n=5000$,  $m=6000$ samples were used as the training set and the validation set has size $1000$.
Regarding REJAFADA, the problem dimension is $n = 6824$. Out of the $1996$ couples $\{(a_i,b_i)\}$, $m = 1597$ couples were used as training set to define problem \eqref{classification} while the remaining $399$ couples were used as validation set.  The corresponding least-squares problem is in this case  underdetermined.
The accuracy in the classification problems is measured as the percentage of labels correctly predicted in the validation set.

\begin{figure}[h]
\centering

        \includegraphics[width = 0.9\textwidth]{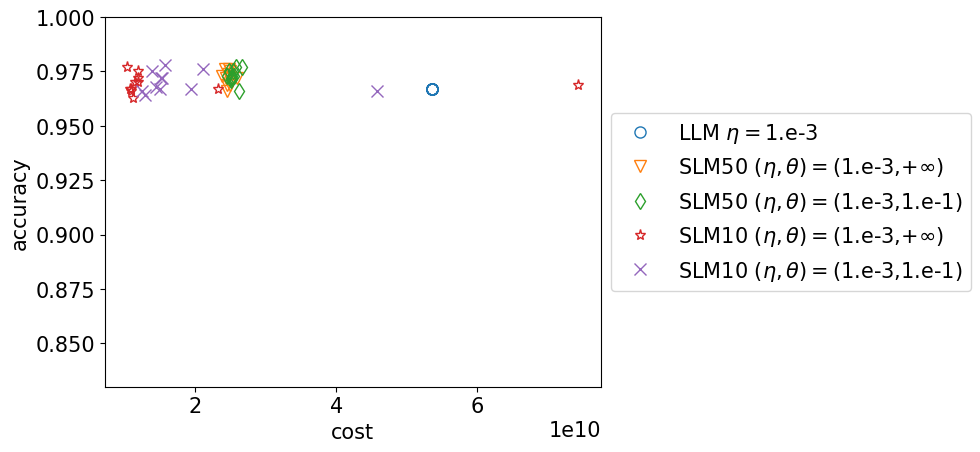}
    \includegraphics[width = 0.9\textwidth]{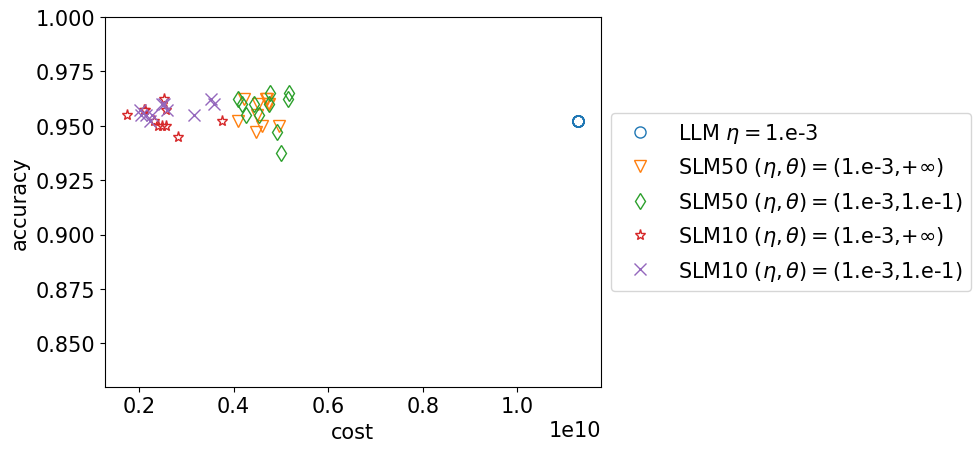}

    \caption{Accuracy at termination versus total computational cost. Upper: GISETTE dataset. Lower: REJAFADA dataset.}
    \label{fig:class1_acc_gradstop_1hashing}
\end{figure}

\begin{figure}[h]
\centering
     \includegraphics[width = 0.9\textwidth]{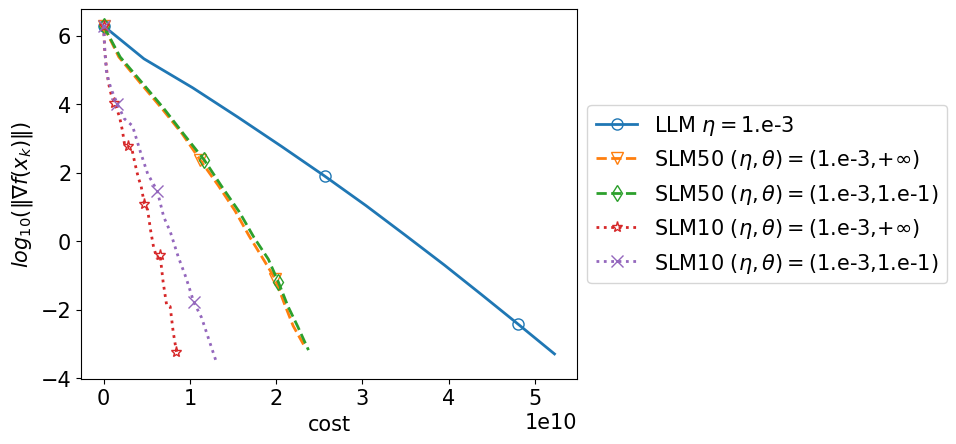} 
     
    \includegraphics[width = 0.9\textwidth]{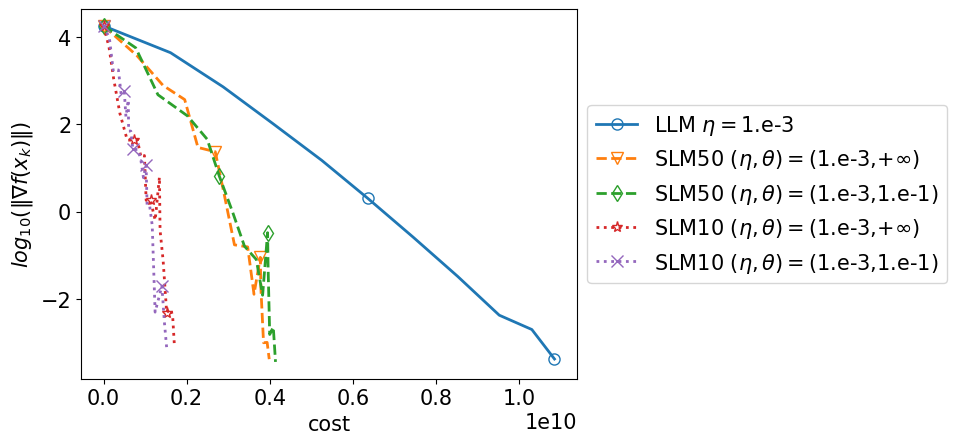}

    \caption{Norm of the gradient $\nabla f(x_k)$ versus computational cost. Upper: GISETTE dataset. Lower: REJAFADA dataset.}
    \label{fig:class1_gradseq_gradstop_1hashing}
\end{figure}

\begin{figure}[h]
\centering

        \includegraphics[width = 0.9\textwidth]{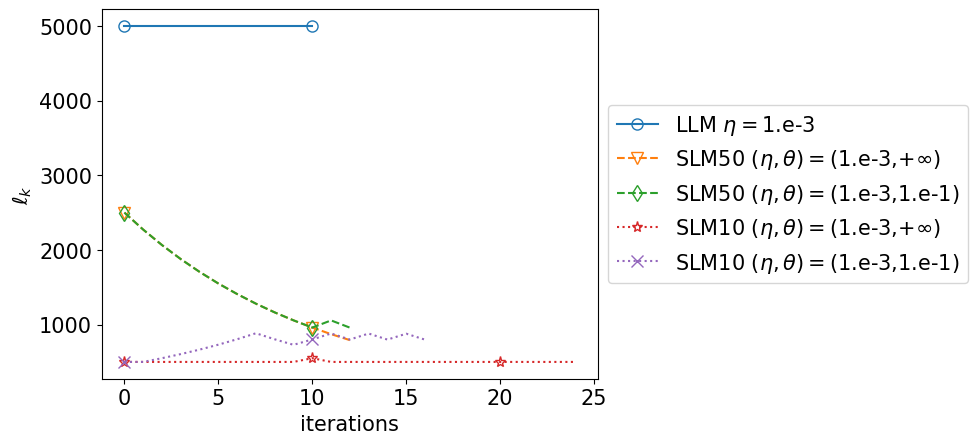}
        
    \includegraphics[width = 0.9\textwidth]{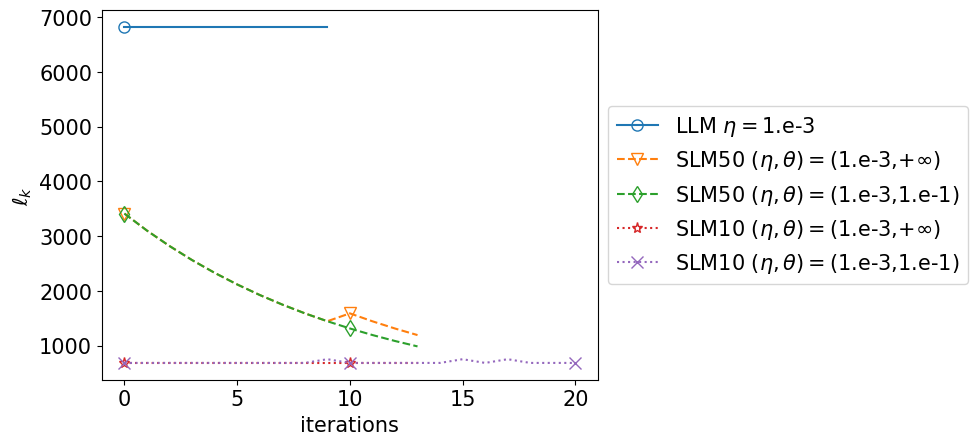}

    \caption{Sketching dimension $\ell_k$ at each iteration. Upper: GISETTE dataset. Lower: REJAFADA dataset.}
    \label{fig:class1_size_gradstop_1hashing}
\end{figure}

We solved this problem with LLM and SLM Algorithms and null initial guess $x_0$. We present results obtained using SLM10 and SLM50 with constant forcing terms  $\eta_k = \eta= 10^{-3}$, $\forall k\ge 0$,  and $\theta \in\{  +\infty, 10^{-1}\}$. 
For every couple of parameters $(\eta, \theta)$, we run SLM  11 times. 
The results obtained are reported in Figures \ref{fig:class1_acc_gradstop_1hashing}--\ref{fig:class1_size_gradstop_1hashing}.

In Figure \ref{fig:class1_acc_gradstop_1hashing}, for each run of LLM and SLM and each pair $(\eta, \theta)$  we plot the accuracy at termination versus  the total computational cost.
The computational cost is defined as follows. Each evaluation of the residual function $F(x_k)$ requires the computation of $m$ scalar products of the form $a_i^T x$ and the overall cost is $mn$. Such scalar products can be stored and used to evaluate the Jacobian $J(x_k)$, whose cost can therefore be disregarded. The evaluation of the gradient $J(x_k)^T F(x_k)$ has cost $mn$ and the evaluation of  $J(x_k)^T J(x_k)s_k$ has cost $2mn$. Finally, each iteration of LSMR method requires two matrix vector products of sizes $m \ell_k$. To summarize, the per-iteration cost is given by $2m\ell_kq_k+4mn$, 
where $q_k$ is the number of LSMR iterations performed at outer-iteration $k$. For the LLM method, such cost is evaluated setting  $\ell_k = n$, $\forall k\ge 0$.   
In figure \ref{fig:class1_gradseq_gradstop_1hashing} we plot the norm of the gradient $\nabla f(x_k)$ versus the computational cost, for the median run in terms of overall computational cost.
Figure \ref{fig:class1_size_gradstop_1hashing} shows how the sketching size $\ell_k$ evolves through the iterations. For each algorithm and pair $(\eta, \theta)$, we plot the results that correspond to the median run in terms of the final computational cost.

In Figure \ref{fig:class1_acc_gradstop_1hashing} we can notice that all runs achieve approximately the same accuracy on the validation set. However, the overall computational cost is significantly smaller for the SLM algorithms compared to the LLM algorithm, except for two runs of SLM10 applied to  GISETTE problem, and the best results are obtained by the SLM10 algorithm.  The savings obtained with the SLM Algorithm are also shown in Figure \ref{fig:class1_gradseq_gradstop_1hashing}.

\vskip 20pt
\appendix
\renewcommand{\appendixname}{}
\section{Appendix}
\subsection{Matrix distributions}\label{Auno}
We summarize the definition of some matrix distributions of interest and specify their parameters with respect to properties \eqref{embedding} and \eqref{bounded_norm}. We denote $M_{i,j}$ the entries of a matrix $M$.

\begin{definition}
    $M\in\R^{\ell\times n}$ is a Scaled Gaussian matrix if its entries $M_{i,j}$ for $i=1, \ldots, \ell,\ j=1,\ldots, n$ are i.i.d. and distributed as $\mathcal{N}(0,\ell^{-1})$.\\
\end{definition}
\begin{definition}\label{def:1&shashing}
  Given $\ell\leq n$, $s\leq \ell$, $M\in\R^{\ell\times n}$ is an $s$-hashing matrix if for every column index $j\in\{1,\dots,n\}$ we sample without replacement  $i_1,\dots i_s$ uniformly at random  and set $M_{i_p,j} = \pm 1/\sqrt{s}$, $p=1,\ldots,s$. \\
\end{definition}

\begin{definition}\label{def::stable-1-hashing}
Let $\ell \in \mathbb{N}^+$ with $\ell<n$. A stable $1$-hashing matrix 
$M \in \R^{\ell \times n}$ has one non-zero per column, whose value is $\pm 1$ with equal probability, with the row indices of the non-zeros given by the sequence $\cal{I}$ constructed as follows.  Repeat the set $\{1,2, \dots, \ell\}$ for $\left\lceil 
{n/\ell}\right\rceil$ times to obtain a set ${\cal{D}}$. Then randomly sample $n$ elements from ${\cal{D}}$ without replacement 
to construct the sequence $\cal{I}$. \\
\end{definition}

\begin{definition}
    $M\in\R^{\ell\times n}$ is a Scaled Sampling matrix if for every row index $i=1, \ldots, \ell$ we sample $j\in\{1,\dots,n\}$ uniformly at random and set $M_{i,j} = \sqrt{n/m}.$\\
\end{definition}

 For the classes of matrices introduced above, Table \ref{tab:matrici} from  \cite[Table 4.2]{Shao21}   summarizes the value of $\varepsilon$, $M_{\max}$, $\delta_M^{(1)}$, $\delta_M^{(2)}$, $\ell$ for the fulfillment of \eqref{embedding} and \eqref{bounded_norm}.
  Notice that \eqref{bounded_norm} holds with probability 1 for all considered matrices except for scaled Gaussian matrices, and that 
  the value  $M_{\max}$ decreases as $\ell$ increases for stable $1$-hashing and scaled sampling matrices.

\begin{table}{\footnotesize
\begin{center}
\begin{tabular}{ c|c|c|c|c|c|c|c} 
 & $\varepsilon$ & $\delta_M^{(1)}$ & $\ell$  & $\delta_M^{(2)}$ & $M_{\max}$\\
 \hline
 Scaled Gaussian & (0,1) & $e^{-\frac{\varepsilon^2\ell}{4}}$  & $\frac{4}{\varepsilon^2}\log\left(\frac{1}{\delta_M^{(1)}}\right)$ & (0,1) &  $1+\sqrt{\frac{n}{\ell}}+\sqrt{\frac{2\log(1/\delta_{M}^{(2)})}{\ell}}$\\
 \hline
 $s$-hashing & (0,1) & $e^{-\frac{\varepsilon^2\ell}{C_1}}$  & $\frac{C_1}{\varepsilon^2}\log\left(\frac{1}{\delta_M^{(1)}}\right)$ & 0 &  $\sqrt{\frac{n}{s}}$
 \\
  \hline
 Stable $1$-hashing & (0,3/4) & $e^{-\frac{(\varepsilon-1/4)^2\ell}{C_1}}$  & $\frac{C_3}{(\varepsilon-1/4)^2}\log\left(\frac{1}{\delta_M^{(1)}}\right)$ & 0 &  $\sqrt{\lceil\frac{n}{\ell}\rceil}$
 \\
  \hline
 Scaled Sampling & (0,1) & $e^{-\frac{\varepsilon^2\ell}{2n\nu^2}}$  & $\frac{2n\nu^2}{\varepsilon^2}\log\left(\frac{1}{\delta_M^{(1)}}\right)$ & 0 &  $\sqrt{\frac{n}{\ell}}$
\end{tabular}
\end{center}
\caption{Values of $\varepsilon, M_{\max}, \delta_M^{(1)}$ and $\delta_M^{(2)}$ in \eqref{embedding} and \eqref{bounded_norm} for different classes of matrices}\label{tab:matrici} 
}
\end{table}

\vskip 5pt
\subsection{Proof from Section \ref{sec3}}\label{Asec3}
\begin{description}
    \item{} 
\begin{proof} 
[Proof of Lemma \ref{lemmaA12cfs}]
$i)$  We proceed by induction on $N$ and consider $N=1$ first.  Since  $= e^{-\lambda x}$ is convex,
we have $e^{-\lambda T_0}\le 1+(e^{-\lambda}-1)T_0$ and 
\begin{equation}
\E\left[e^{-\lambda T_0} \right]\le 1+(e^{-\lambda}-1)\E\left[ T_0\right].
\label{eqn::inductionZeroUpper}
\end{equation}
Moreover, we have        $ \E \left[T_0 \right]\geq \mathbb{P}(T_0=1) \geq 1-\delta_M$,
where the first inequality is due to $T_0\geq 0$, and the second inequality to Assumption \ref{Atrue}.
Therefore, noting that $e^{-\lambda}-1<0$, \eqref{eqn::inductionZeroUpper} gives
\begin{equation}\label{eqn::inductionFirstConclusion}
\E \left[e^{-\lambda T_0} \right]
\leq 1 + (e^{-\lambda}-1) (1-\delta_M) \leq e^{(e^{-\lambda}-1) (1-\delta_M)}, 
\end{equation}
where the last inequality comes from $1+y \leq e^y$ for $y \in \R$.  

Now assume 
$
\E \left[e^{-\lambda \sum_{k=0}^{N-2} T_k } \right] \leq 
\left[ e^{( e^{-\lambda}-1)(1-\delta_M)}\right]^{N-1}  .  
$
Due to the Tower property, we have 
\begin{align}
& \E \left[e^{-\lambda \sum_{k=0}^{N-1} T_k}\right]= \E \left[ \E\left[e^{-\lambda \sum_{k=0}^{N-1} T_k}\,|\, {\cal{F}}_{N-2}\right] \right],
\notag 
\end{align}
and 
\begin{align}
\E\left[e^{-\lambda \sum_{k=0}^{N-1} T_k}| {\cal{F}}_{N-2}\right]&=  \E\left[\prod_{k=0}^{N-1} e^{-\lambda T_k}\,|\, {\cal{F}}_{N-2}\right]  \notag \\
&= \prod_{k=0}^{N-2} e^{-\lambda T_k}\E\left[  e^{-\lambda T_{N-1}}\,|\, {\cal{F}}_{N-2}\right]\\
&\le e^{( e^{-\lambda}-1)(1-\delta_M)} \prod_{k=0}^{N-2} e^{-\lambda T_k},
\end{align}
since $T_{N-1}$ is conditionally independent of the past iterations $T_0, \ldots, T_{N-1}$, and 
in the last inequality we used (\ref{eqn::inductionFirstConclusion}) and the arguments for the case $N>1$
(from
Assumption \ref{Atrue}, $\E \left[T_{N-1}\,|\, {\cal{F}}_{N-2} \right] =\mathbb{P}(T_{N-1}=1\,|\, {\cal{F}}_{N-2}) \geq 1-\delta_M$).
Hence, induction implies
\begin{eqnarray*}
\E \left[e^{-\lambda \sum_{k=0}^{N-1} T_k}\right] &\le& 
e^{( e^{-\lambda}-1)(1-\delta_M)} \E\left[\prod_{k=0}^{N-2} e^{-\lambda T_k}\right]\\
&\le& e^{( e^{-\lambda}-1)(1-\delta_M)} \left[ e^{( e^{-\lambda}-1)(1-\delta_M)}\right]^{N-1} ,   
\end{eqnarray*}
and the claim in Item $i)$ follows.
\vskip 5pt \noindent
$ii)$ See \cite[Proof of Lemma A.1]{cfs}.
\end{proof}
\end{description}
\vskip 5pt
\subsection{Proofs from Section \ref{sec5}}\label{Asec5}
\begin{description}
    \item{} 
\begin{proof} 
[Proof of Lemma \ref{lemma:loc_s}] By \eqref{bound_sv} and Assumption \ref{ass:eigJ} we have 
 $\lambda^{r_k}_k\geq (1-\varepsilon)\lambda_r(J(x_k)^T J(x_k))\geq (1-\varepsilon)\lambda_{\min}.$ Using \eqref{bound_s_hat}, the fact that iteration $k$ is true,  Item  4 in Lemma \ref{lemma_ineq}, and the assumption on $\eta_k/\mu_k$ we have
 \begin{equation}\begin{aligned}
     \|s_k\|_2 =& \|M_k^T\widehat{s}_k\|_2\leq M_{\max}\left(\frac{1}{(1-\varepsilon)\lambda_{\min}+\mu_k} +\frac{\eta_k}{\mu_k}\right)\|M_kJ(x_k)^T F(x_k)\|_2\\&\leq M_{\max}(1+\varepsilon)^{1/2}\left(\frac{1}{(1-\varepsilon)\lambda_{\min}+\mu_k} +\frac{\eta_k}{\mu_k}\right)\|J(x_k)^T F(x_k)\|_2\\
     &\leq M_{\max}(1+\varepsilon)^{1/2}\left(\frac{1}{(1-\varepsilon)\lambda_{\min}} +\bar c\right)L \dist(x_k,\Omega^*),
 \end{aligned}\end{equation}
 therefore the thesis holds with 
  $c_1 = M_{\max}(1+\varepsilon)^{1/2}\left(\frac{1}{(1-\varepsilon)\lambda_{\min}} + \bar c\right)L .$
\end{proof}
\item{}

\vskip 5pt
\begin{proof} [Proof of Lemma \ref{lemma:loc_dist}]
    By Assumption \ref{ass:errorbound}, the definition of $\theta_k^*$ in (\ref{thetak}), and Item 4 in Lemma \ref{lemma_ineq}, we have
    \begin{equation}
    \begin{aligned}                \omega\dist(x_{k+1},\Omega^*)&\leq \|\nabla f(x_{k+1})\|_2\leq \|\nabla f(x_{k+1})-\nabla f(x_k) - J(x_k)^T J(x_k) s_k\|_2 + \\
    & \quad \,  \, \|\nabla f(x_k) +  J(x_k)^T J(x_k) s_k\|_2  \\
    &=\|\nabla f(x_{k+1})-\nabla f(x_k) - J(x_k)^T J(x_k) s_k\|_2 + \theta^*_k \|J(x_k)^TF(x_k)\|_2 \\ 
    &\leq \|\nabla f(x_{k+1})-\nabla f(x_k) - J(x_k)^T J(x_k) s_k\|_2 + \theta^*_k L \dist(x_k,\Omega^*).
    \end{aligned}
    \end{equation}
    Using  Items 2 and 3  in Lemma \ref{lemma_ineq},  \eqref{bound_ass5.4}
    and proceeding as in Lemma 4.1 in \cite{santos}, we get the thesis, with
    $L_3 = \sigma(1+(1+c_1)^{1+\beta})$ and
    $L_4 = L_2c_1^2+L_0{J}_{\max}(1+c_1)(2+c_1)$.
\end{proof}
\end{description}

\end{document}